\documentclass[aip,jcp,preprint]{revtex4-1}

\usepackage[width=16cm,  top=2.5cm, footskip=2cm]{geometry}

\usepackage{amsfonts}
\usepackage{amssymb}       
\usepackage{amsmath}
\usepackage{bm}

\usepackage{lineno}

\usepackage[dvips]{graphicx}
\usepackage{graphicx}

\usepackage{booktabs}
\usepackage{hyperref}
\usepackage{breakurl}

\usepackage{verbatim}
\usepackage{subeqnarray}

\usepackage{color}

\usepackage[normalem]{ulem}

\def\LOPER{\mathcal{L}}
\def\LOPERa{\mathcal{L}_A}
\def\LOPERb{\mathcal{L}_B}
\def\LOPERcpl{\mathcal{\bar L}}
\def\LOPERcplSimple{\mathcal{L}^0}
\def\LATT{{\Lambda}_N}
\def\RATEa{c_A}
\def\RATEb{c_B}
\def\MEANUa{u_A}
\def\MEANUb{u_B}
\def\PROCa{\sigma_t}
\def\PROCb{\eta_t}
\def\PROCMICa{\{\sigma_t\}_{t\geq 0}}
\def\PROCMICb{\{\eta_t\}_{t\geq 0}}
\def\PROCa{\sigma_t}
\def\PROCb{\eta_t}
\def\R{\mathbb{R}}
\def\FISHERR{\mathbf{F_{\mathcal{H}}}}

\def\RELENT#1#2{\mathcal{R}\left({#1}\SEP{#2}\right)}
\def\SEP{{\,|\,}}           

\def\SPINSP{\Sigma}                     
\def\CONFNEW{\sigma^{x,\omega}}
\def\SEP{{\,|\,}}           
\def\PERIOD{\,.}            
\def\SIGMA{{\mathcal{S}}}               
\def\NEIGH{N}	
\def\CONF{\Omega}
\def\NUMSP{{K}}                       
\def\VIZ#1{(\ref{#1})}      
\def\PROB#1{{\mathbb{P}\left({#1}\right)}}

\def\EXPECT{{\mathbb{E}}}
\def\VAR{\operatorname{Var}}

\def\SGN{\operatorname{sign}}

\def\EXP#1{e^{#1}}
\def\PROC#1{\{{#1}\}_{t\geq 0}}
\def\BIGO{\mathcal{O}}

\def\COMMA{\,,}             

\def\NORMB#1{\Vert\,#1\,\Vert_\infty}

\newtheorem{definition}{Definition}[section]
\newtheorem{remark}{Remark}[section]
\newtheorem{optProb}{Optimization Problem}

\begin{document}

\title{Goal-oriented sensitivity analysis for lattice kinetic Monte Carlo simulations}

\author{Georgios Arampatzis}
\affiliation{Department of Applied Mathematics, University of Crete and  Foundation of Research and Technology-Hellas, Greece,
{\tt garab@math.uoc.gr}}

\author{Markos A. Katsoulakis}
\affiliation{Department of Mathematics and Statistics, University of Massachusetts, Amherst, MA, USA,
{\tt markos@math.umass.edu}}

\begin{abstract}
In this paper we propose a new class of coupling methods for the sensitivity analysis of high dimensional stochastic systems and in particular for lattice Kinetic Monte Carlo. 
Sensitivity analysis for stochastic systems  is typically based on approximating continuous derivatives with respect to model parameters by the mean value of samples from a finite difference scheme. Instead of using independent samples the proposed algorithm reduces the variance of the estimator by developing a   strongly correlated-"coupled"- stochastic process for both the perturbed and unperturbed stochastic processes,  
defined in a common state space. The novelty of our construction is  that the  new coupled process depends on the 
 targeted observables, e.g. coverage, Hamiltonian, spatial correlations, surface roughness, etc., 
hence we refer to the proposed method   as  {\em goal-oriented} sensitivity analysis.
In particular, the rates of the coupled  Continuous Time Markov Chain  are obtained as  solutions to  a goal-oriented  optimization problem, depending  on the observable of interest,
 by considering the minimization functional of  the corresponding variance. We show that this functional can be  used as a diagnostic tool for the design and evaluation of different classes of couplings. Furthermore, the resulting KMC sensitivity algorithm has an easy implementation that is based on the Bortz--Kalos--Lebowitz algorithm's philosophy, where here events are divided in classes depending on level sets of the observable of interest. Finally,  we demonstrate in 
several examples including adsorption, desorption and diffusion Kinetic Monte Carlo that  for the same confidence interval and observable, the proposed  goal-oriented algorithm can be  two orders of magnitude faster than existing  coupling algorithms for spatial KMC such as the  Common Random Number approach.
\end{abstract}

\maketitle


\section{Introduction}\label{section_introduction}

Recently there has been significant progress in developing sensitivity analysis
tools for stochastic processes
modeling well-mixed chemical reactions and biological networks. Some of the mathematical
tools include  log-likelihood methods and Girsanov transformations \cite{Glynn:90,Nakayama:94,Plyasunov:07}, 
polynomial chaos \cite{Kim:07},   finite difference methods and their variants 
 \cite{Rathinam:10, Anderson} and pathwise
sensitivity methods \cite{Khammash:12}; a somewhat parallel literature, facing related challenges, 
exists also in mathematical finance \cite{Glassermann} and operations research \cite{CaoChen,Dai,Bremaud}.  However, existing
sensitivity analysis approaches    can have 
an overwhelming computational cost in high dimensions, such as lattice Kinetic Monte Carlo (KMC), either due to high variance in the gradient estimators,
or in models with a high-dimensional parameter space.
Such issues and  comparisons between methods are 
discussed, for instance, in recent literature \cite{Vlachos:12, Srivastava:13, PK_RER}, see also the demonstration in Figure~\ref{fig:simulation1_variance}.

Estimating the sensitivity of a stochastic  process $\{\sigma_t=\sigma_t(\theta),\, t\ge 0\}$ with respect to perturbations in the model parameters $\theta=(\theta_1,..., \theta_k)$ can be viewed
either at the level of the process'  probability distribution function (PDF) %
or as a response 
of specific averaged observables $f=f(\sigma_t)$ in some time interval $t \in [0, T]$, namely  the quantity
$$u(t,\rho;\theta)=\EXPECT_\rho[f(\sigma_t(\theta))]\COMMA$$
to parameter perturbations, where $\rho$ denotes the initial configuration of the stochastic  process $\sigma_t$ and $\EXPECT_\rho$ the corresponding expected value. We focus on the latter  perspective, which we also refer to 
as a {\it goal-oriented} approach, since the focus is on observables of interest.
In this case, we  quantify the sensitivity of specific observables %
by
estimating gradients of the type $\partial_{\theta_l} u(t,\rho;\theta)$, where $l \in \{1, ..., k\}$ and $\theta=(\theta_1,..., \theta_k)$.
%
%
In turn, it is   commonplace to evaluate  such gradients for  parametric sensitivity of the observable 
$u(t,\rho;\theta)$
 by using finite difference (FD) approximations, e.g., for   first order derivatives
\begin{equation}\label{partial_derivative}
\partial_{\theta_l} u(t,\rho;\theta) 
\approx \Delta := \frac{1}{h} (u(t,\rho;\theta + h e_l) - u(t,\rho;\theta))  \COMMA
\end{equation}
where $h\in\mathbb{R}$ is a small parameter and $e_l$ a unit vector in $\mathbb{R}^k$ with $e_{l,j}=1$ if $l=j$ and zero otherwise. 
While such simulation approach appears straightforward,  it suffers from a well-known problem arising from the high variance of the estimator 
for the finite-difference $\Delta$ \cite{Rathinam:10}, at least if we naively pick independent samples as we further explain in (\ref{joint_probability}) below. 

Overall, developing  methods of reduced variance is a {\em critical computational task} for carrying out sensitivity analysis in high-dimensional
complex stochastic models such as the spatial Kinetic Monte Carlo algorithms studied here.
The high computational cost of individual realizations of the stochastic process renders prohibitive the generation of a large number of samples for reliable  ensemble averaging.
Hence, reducing the variance by  orders of magnitude, e.g.  Figure~\ref{fig:simulation1_variance}, will result in an analogous 
reduction in the number of samples needed for the ensemble average.%
%

%
%
%
%
The variance of finite difference  sensitivity  estimators is usually reduced by  employing {\em coupling methods},
i.e., by constructing highly correlated paths for the processes  $\sigma_t(\theta)$ and $\sigma_t(\theta +h e_l)$.
The simplest
coupling is to run the two processes using the same stream of random number, known as Common Random Number (CRN) method,
which for spatially extended systems the induced correlation is not enough to reduce the variance.
Improvements of this method is the Common Random Path algorithm \cite{Rathinam:10} (CRP)   and the Coupled Finite Difference (CFD) method \cite{Anderson}.
In fact,
 it is shown  \cite{Srivastava:13} that among these  coupling methods CFD  performs better, at least for relatively simple, low-dimensional reaction networks.  
In general, such    couplings  are  suitable only for low-dimensional, well-mixed systems as  variance estimates depend on system size \cite{Anderson} and more importantly  
it was not clear, up to now,  how to extend them in an efficient manner to   spatially distributed models Kinetic Monte Carlo (KMC) models.
As we see in  Figure~\ref{fig:simulation1_variance} depicting a lattice KMC simulation of a spatially extended adsorption/desorption model, 
the variance of such coupled  estimators remains very high even if CRN is implemented.
%
%

%
%
The novelty of the approach we propose  in this direction relies on developing  a different concept of stochastic coupling which (a) is
 suitable for spatially extended  systems such as lattice KMC,  and  (b) is  designed for specific observables $\EXPECT_\rho[f(\sigma_t)]$ hence  we call it a goal-oriented coupling method.
%
%
Our proposed method relies  on defining a new coupled continuous time Markov Chain through a suitable generator that acts on observables of the involved high-dimensional stochastic processes associated with parameters $\theta$ and $\theta+\epsilon$, where $\epsilon$ denotes any $k$-dimensional perturbation in parameter space, namely
$$
\epsilon = h e_l, \quad l=1,\ldots,k \PERIOD
$$
Indeed, we define the  coupling  $\zeta_t = (\PROCa,\PROCb)$ of two    stochastic processes $\sigma_t = \sigma_t(\theta)$  and  $\eta_t=\sigma_t(\theta+\epsilon)$, i.e.  we couple the dynamics, setting them    in a common, product  probability space,   but at the same time  we  respect 
the marginal distributions  of each one of them. 
Clearly we have freedom on how to select these coupling rates, however our goal is to minimize the variance of estimators of the finite difference of specific observables $f$, while  keeping the computational cost of coupling low. 
The rates of the coupled stochastic processes $\zeta_t = (\PROCa,\PROCb)$ are obtained by solving an optimization problem associated with  minimizing  the  variance between the coupled stochastic dynamics. 
Furthermore, the   optimization functional is also a measure of the ``tightness'' of the coupling, allowing us to compare various coupling schemes and systematically assess their effectiveness in variance reduction.
The algorithmic implementation of the proposed coupled method is a Bortz-Kalos-Lebowitz (BKL)--type algorithm  in the sense that events are categorized into pre-defined sets. In the classical BKL algorithm \cite{BKL}  events are divided into classes  of equal rates, i.e.  according to level sets of the rates. However, here the events are divided into classes depending on the  observable's level sets, since we are interested in tight couplings of time series of specific observables. 
%
Numerical examples of spatial KMC, e.g. adsorption/desorption/diffusion processes, are presented throughout  the paper and demonstrate   that the variance  can be improved by  two orders of magnitude compared to coupling methods used up to now, such as the CRN method. At the same time the computational overhead of computing the coupled rates is two times slower than that of the CRN method leading to an overall speed up factor of two orders of magnitude. 
Furthermore, numerical experiments also demonstrate that the variance-related optimization functional indeed
constitutes a diagnostic tool for the design and evaluation of different couplings.

The paper is structured as follows: in Section 2 we provide  background and notation for spatial KMC methods and discuss earlier work on coupling methods. 
In Section 3 we introduce coupling methods for spatial KMC algorithms
and demonstrate the resulting variance reduction  in several examples. In Section 4 we introduce  improved coupling algorithms which do not attempt to couple the time series configurations of the entire stochastic process, but instead they are designed to couple only time series of  specific observables; the resulting algorithms are  constructed in the  spirit of the BKL algorithms for KMC and we demonstrate variance reduction up to two orders of magnitude. In Section 5, 
we discuss the limitations, as well as the potential applicability of the proposed coupling  to  systems with a very high-dimensional parameter space. 
Finally, in Appendix  we give detailed examples of complex reaction-diffusion models and a description of the implementation of the proposed sensitivity analysis method.


\section{ Background }

\subsection{ Markov Chains and kinetic Monte Carlo}\label{section:background_MC}

In this work we present the proposed sensitivity analysis methods in the context 
of spatial Kinetic Monte Carlo (KMC) methods,  although similar challenges and ideas are relevant 
to all other molecular simulation methods.
%
%
%
%
%
%
%
The resulting stochastic processes studied in this work are set on a discrete, albeit high-dimensional, configuration space  $\SIGMA$ and necessarily
have to be of jump type describing transitions between different configurations $\sigma \in \SIGMA$. Mathematically, such a Continuous Time Markov Chain (CTMC)  is a stochastic
process $\{\sigma_t\}_{t\geq 0}$ defined completely in terms of  the local  transition rates $c(\sigma, \sigma')$ which determine  the 
updates (jumps)  from any current state $\sigma_t=\sigma$ to a (random) new state $\sigma'$. Usually, simulated paths of the process are
constructed via Kinetic Monte Carlo  (KMC), that is  through the procedure described in \VIZ{totalrate} and \VIZ{skeleton} below. 

Realizations of the process are constructed from the embedded discrete time  Markov chain \cite{KL} $S_n =
\sigma_{t_n}$  with jump times $t_n$ from the exponential distribution: 
In the context of the spatially distributed problems  (in which we are interested  here), the  local transition rates will be denoted as
$c(\sigma, \sigma';\theta)$ where $\theta\in\R^k$ is a vector of the model parameters, describing transitions from
the configuration at time $t$,  $\sigma_t=\sigma$ into a new configuration $\sigma'$ .
The local transition rates $c$, define the total rate
\begin{equation}\label{totalrate}
\lambda(\sigma;\theta)=\sum_{\sigma'} c(\sigma,\sigma';\theta)\COMMA
\end{equation}
which is the intensity of the exponential waiting time for a jump to be performed when the system is currently at the state $\sigma$. 
The transition probabilities for the embedded Markov chain $\{S_n\}_{n\geq 0}$ are 
\begin{equation}\label{skeleton}
   p(\sigma, \sigma';\theta)=\frac{c(\sigma,\sigma';\theta)}{\lambda(\sigma;\theta)}\PERIOD
\end{equation}
In other words once the exponential ``clock'' signals a jump, the system transitions from the state $\sigma$ 
to a new configuration $\sigma'$ with probability $p(\sigma, \sigma';\theta)$.
On the other hand, the evolution of the entire system at any time $t$ is described by the  transition probabilities
$P(\sigma, t ; \sigma_0;\theta):=\PROB{\sigma_{t} = \sigma\SEP \sigma_0 = \sigma_0}$
where $\sigma_0 \in \SIGMA$ is any  initial configuration. The transition probabilities, corresponding to the local rates 
$c$, satisfy the Forward Kolmogorov Equation \cite{Gardiner04} ({Master Equation}) ,
\begin{equation}\label{master}
\partial_t P(\sigma, t ; \sigma_0;\theta):=\sum_{ \sigma'\ne \sigma} c(\sigma', \sigma;\theta)P(\sigma', t; \sigma_0;\theta)-c(\sigma,\sigma';\theta)P(\sigma, t ;
\sigma_0;\theta) \COMMA
\end{equation}
where $P(\sigma, 0 ; \sigma_0;\theta)=\delta (\sigma-\sigma_0)$ and $\delta (\sigma-\sigma_0)=1$ if $\sigma=\sigma_0$
and zero otherwise.

\noindent {\bf Generators for CTMC.}   Typically in KMC we need to compute expected values of such observables, that is quantities defined as
\begin{equation}
   u(\sigma_0, t):=\EXPECT_{\sigma_0}[f(\sigma_t)]=\sum_{\sigma'} f(\sigma') P(\sigma', t; \sigma_0)\, \PERIOD
\label{propagator1}
\end{equation}
Here $\EXPECT_{{\sigma_0}}$ denotes the expected value with respect to the law 
of the process $\{\sigma_t\}$ conditioned
on the initial configuration ${\sigma_0}$. 
By a straightforward calculation \cite{Gardiner04} using \VIZ{master}  we obtain that  the observable \VIZ{propagator1}
satisfies the initial value problem
\begin{equation}\label{ODE}
\partial_t u({\sigma_0}, t)=\LOPER u({\sigma_0}, t)\, , \quad\quad u({\sigma_0}, 0)=f({\sigma_0})\, ,
\end{equation}
where the operator $\LOPER$ is known as the {\em generator} of the CTMC \cite{Liggett}
\begin{equation}\label{generator}
 \LOPER f(\sigma) = \sum_{\sigma'}c(\sigma, \sigma')[f(\sigma')-f(\sigma)] \PERIOD
\end{equation}
The generator fully determines the process $\sigma_t$ while (\ref{generator}) can viewed as the dual \cite{Liggett} of (\ref{master}).
Although in order to describe KMC algorithms (and in general any discrete space continuous time Markov process) is 
not necessary to use generators, we have found that here they allow us to systematically construct and assess couplings 
of stochastic processes for the purpose of providing low variance finite difference estimators. For this reason we introduce 
the concept, hopefully in a self contained manner, and exploit it in Sections \ref{sec:microscopic_coupling} and \ref{sec:macroscopic_coupling}.
Using generator notation \cite{Liggett} we then can rewrite \VIZ{propagator1}, as the the action of the Markov semi-group  propagator $\EXP{t\LOPER}$
associated with the generator $\LOPER$ and the process $\PROC{\sigma_t}$ acting
on the observable $f$:
\begin{equation}\label{semigroup}
    u({\sigma_0}, t)=\EXPECT_{\sigma_0}[f(\sigma_t)] := \EXP{t\LOPER}f({\sigma_0})\, \PERIOD
\end{equation}
In fact, in the case where the state space is finite but high dimensional, as in spatial KMC processes considered here, the operator $\LOPER$ is essentially a matrix (bounded), that is
\begin{equation}\label{bounded}
 \NORMB{\LOPER f} \leq c \NORMB{f}
\end{equation}
where $\NORMB{f} = \max_{\sigma\in\SIGMA} |f(\sigma)|$ and $c$ is a constant independent of $f$. Due to the boundness of $\LOPER$, $\EXP{t\LOPER}$ can be also defined as an infinite series.  For $t=\delta t \ll 1$ the
semigroup \VIZ{semigroup} can be approximated by series truncation as 
\begin{equation}\label{approximation_semigroup}
e^{ \delta t \LOPER} = I + \delta t \LOPER + \BIGO (\delta t^2) \COMMA
\end{equation}
allowing us to write the solution \cite{KL} of \VIZ{ODE} in terms of the
operator $\LOPER$ and $\delta t$,
\begin{equation}\label{solution_series}
u( \delta t;\sigma)  =  u(0 ;\sigma) + \delta t \LOPER u(0 ;\sigma) + \BIGO (\delta t^2) \PERIOD
\end{equation}
Rigorous statements with less stringent conditions on the types of such processes and corresponding generators can be found in literature \cite{Liggett}. 

\subsection{KMC on a Lattice and Benchmark Examples	}\label{sec:lattice_kmc}

We consider interacting particle system defined on a $d$-dimensional lattice $\LATT$ of any type (square, hexagonal, Bravais,  etc), where $N$ is the size of the lattice.
As a result, we model the dynamics of the configuration space with a Continuous Time Markov Chain jump process.
We restrict our discussion to lattice gas models where the order parameter or the spin variable
takes values in a finite set $\SPINSP=\{0,1,\dots, \NUMSP \}$.
At each lattice site $x\in \LATT$ an order parameter (or a spin variable in Ising systems) $\sigma(x)\in \SPINSP$ 
is defined. The elements in $\SPINSP$ correspond to occupation of the site $x\in\LATT$ by different
species. The stochastic process $\PROCMICa$  takes values in the  configuration space  
\begin{equation*}
 \SIGMA = \Bigl\{ \bigr(\sigma(x_1),\ldots,\sigma(x_N)\bigl) \;|\; x_i\in\Sigma, \, i=1,\ldots,N \Bigr \} \PERIOD
\end{equation*}
Microscopic dynamics are described by changes (transitions) of spin variables at different sites . 
We study systems in which the transitions are localized and involve only a finite number of sites at each
transition step. First, the {\em local} dynamics are described by an updating mechanism and corresponding  transition rates
$c(x,\sigma)$, such that 
the configuration at time $t$,  $\sigma_t=\sigma$ changes into a new configuration $\sigma'=\sigma^x$ by an
update of $\sigma$ at  the site $x\in\LATT$.
For example, if $\SPINSP=\{0,1\}$ the order parameter models the classical lattice gas with
a single species occupying the site $x$ when $\sigma(x)=1$ and with the site being vacant if $\sigma(x)=0$. The system can jump from state
$\sigma$ to the new configuration $\sigma^x$ (adsorption $0\rightarrow 1$, desorption $1\rightarrow 0$), where\cite{KV_2003,Liggett}
\begin{equation}\label{ising_new_configuration}
 \sigma^x (y) = 
 \begin{cases}
  1-\sigma(x), & \textrm{if } x=y \\
  \sigma(y),   & \textrm{if } x\neq y \COMMA
 \end{cases} 
\end{equation}
and the generator of the process is,
\begin{equation}\label{ising_generator}
   \LOPER f(\sigma) = \sum_{x\in\LATT} c(x,\sigma) \bigl [f(\sigma^x)-f(\sigma) \bigr ] \PERIOD
\end{equation}
The Ising model, with adsorption/desorption dynamics and spins in $\{0,1\}$ has rates,
\begin{equation}\label{ising_rates}
 c_I(x,\sigma) = 
 \begin{cases}
 	c_a, & \textrm{if } \sigma(x) = 0 \\
    c_d e^{-\beta \bigl ( J \left ( \sigma(x-1)+\sigma(x+1) \right) - h \bigr )}, & \textrm{if } \sigma(x) = 1
  \end{cases}
\end{equation}
where $\beta, J$ and $h$ are the inverse temperature, inter-particle potential and external field respectively and $c_a,c_d$ are the adsorption and desorption constants respectively.
Using the same setup as in the Ising model, the simple diffusion process can be described by the transition from the configuration $\sigma$ the new configuration state $\sigma'=\sigma^{x,y}$, where the particle exchanges position with an empty site,
\begin{equation}\label{diffusion_new_configuration}
 \sigma^{x,y}(z) = 
 \begin{cases}
  \sigma(x), & \textrm{if } z=y \\
  \sigma(y), & \textrm{if } z=x \\
  \sigma(z),   & \textrm{if } z\neq x,y \COMMA
 \end{cases} 
\end{equation}
with rate,
\begin{equation}\label{diffusion_rates}
 c_D(x,y,\sigma) = 
 \begin{cases}
 	c_{diff}, & \textrm{if } |x-y|=1 \\
 	0,         & \textrm{if } |x-y| \neq 1 \PERIOD
  \end{cases}
\end{equation}
Modelling interacting diffusions is also straightforward \cite{refmat8}.
Combining simple mechanisms one can build more complex models.
For example the generator that describes an adsorption/desorption type model with diffusion is given by the sum of the generator of the two basic processes,
\begin{equation}\label{ising_diffusion_generator}
\LOPER f(\sigma) = \sum_{x\in\LATT} c_I(x,\sigma) \bigl [f(\sigma^x)-f(\sigma) \bigr ] + \sum_{x\in\LATT}\sum_{y\in\LATT} c_D(x,y,\sigma) \bigl [f(\sigma^{x,y	})-f(\sigma) \bigr ] \PERIOD3
\end{equation}
A  general formulation for the description of more complex  lattice systems and mechanisms that include all of the above is presented in Section
\ref{section:general_form_of_the_generator}, see also Appendix \ref{app:examples}.
These  models describe  a wide range of applications ranging from crystal growth, to catalysis, to biology \cite{refmat8}.
For instance, in catalysis, KMC describes microscopic events, such as adsorption of species on a catalyst and its reverse (desorption to the fluid), 
surface diffusion, and surface reactions; typically these mechanisms occur concurrently. 
In such KMC models the transition rate constants were determined in 
an ad hoc or a semi-empirical manner. However, more recently 
{\em  first-principles} KMC methods  were developed, where kinetic parameters of micro mechanisms  
are estimated by {\em ab initio} density functional theory \cite{Kohn} (DFT) see for instance the recent works 
\cite{HansenNeurock99, Reuter1, evans09, Wu12}. 
Contrary to earlier work, more qualitative descriptions, such first-principles models yield 
a remarkable agreement with experiments \cite{SV2012, Metiu, Christensen}. 
In view of the significant role of DFT based parameter fitting the role of sensitivity analysis is a crucial step in building reliable predictive KMC algorithms.


\subsection{Finite Difference sensitivity and Coupling Methods}
We first consider the following family or generators parametrized by the  parameter vector $\theta \in \R^k$,
\begin{equation}\label{generator_1}
  \LOPER^\theta f(\sigma) = \sum_{x\in\LATT} c(x,\sigma;\theta) \bigl [f(\sigma^x)-f(\sigma) \bigr ] \COMMA
\end{equation}
while we address general lattice KMC in section \ref{section:general_form_of_the_generator}. Our goal is to assess the sensitivity of $u^\theta$, i.e.\ the solution of  \VIZ{ODE} with generator $\LOPER^\theta$, in perturbations in
the parameter $\theta$ of the rates function $c$. Therefore we consider the quantity,
\begin{equation}\label{difference}
 D_{\epsilon}(t;\sigma,\eta)  = u^\theta(t,\sigma) - u^{\theta+\epsilon}(t,\eta) \PERIOD
\end{equation}
%
%
%
From now on, in order to keep the generality and lighten the notation, we define $\MEANUa:=u^{\theta}$ and $\MEANUb:=u^{\theta+\epsilon}$.
The subscript $A$ and $B$ will also be used in rates and generators to correspond to $\theta$ and $\theta+\epsilon$.
In order to estimate the quantity (\ref{difference}), one must simulate many realizations of  the stochastic model, $\sigma^{[i]}_t$,  and
then take ensemble averages to compute $\MEANUa$ and $\MEANUb$. This can be written as
\begin{equation}\label{estimator_mean}
D_{\epsilon}(t;\sigma,\eta) \approx  D_{N,\epsilon}(t;\sigma,\eta) = \sum_{i=1}^N  \frac{ f(\PROCa^{[i]} ) - f( \PROCb^{[i]} )  }{N},
\end{equation}
where $\PROCa^{[i]}$ and $\PROCb^{[i]}$ are the $i-$th  path  at time $t$, generated from $\LOPERa$ and $\LOPERb$
respectively. 

Next we focus on calculating the variance of estimator (\ref{estimator_mean}). First, we write $D_{\epsilon}$ in terms of the probability distributions as,
\begin{equation}
\begin{split}
 D_{\epsilon}(t;\sigma,\eta) &= u^A(t,\sigma_0) - u^B(t,\eta_0) \\ 
 &= \EXPECT_{\sigma_0} f(\sigma_t) -  \EXPECT_{\eta_0} f(\eta_t) \nonumber \\
 &= \sum_{\sigma'} f(\sigma') p^{A}(\sigma',t;\sigma_0)  - \sum_{\eta'} f(\eta') p^{B}(\eta',t;\eta_0) \nonumber \\
 &= \sum_{\sigma',\eta'} \Bigl [ f(\sigma') - f(\eta') \Bigr ] \bar p(\sigma',\eta',t;\sigma_0,\eta_0) \COMMA
\end{split}
\end{equation}
where $p^{A}(\cdot,t;\sigma_0)$ and $p^{B}(\cdot,t;\eta_0)$ are the transition probability measures at time $t$ given that the initial
state is $\sigma_0$ and $\eta_0$, according to  the dynamics of the rate functions $c^A$ and $c^B$ respectively. Here the summation is considered over the entire configuration space, while $\bar p$ denotes the joint probability of the processes  $\PROCMICa$ and $\PROCMICb$ defined as
\begin{equation}
\bar p(\sigma,\eta,t;\sigma_0,\eta_0) = \mathbb{P} (\sigma_t = \sigma,\eta_t=\eta | \sigma_0=\sigma_0,\eta_0=\eta_0) \PERIOD
\end{equation}
Thus, the variance of estimator (\ref{estimator_mean}) is,
\begin{equation}\label{variance_int}
\begin{split}
0  &\leq  \VAR_{(\sigma_0,\eta_0)} \bigl [  D_{N,\epsilon}(t;\sigma,\eta)  \bigr ] =  \VAR_{(\sigma_0,\eta_0)} \bigl [ f(\PROCa) -  f(\PROCb) \bigr ] \\
&= \sum_{\sigma',\eta'} \Bigl [ \bigl(f(\sigma')-\EXPECT_{\sigma_0} f(\sigma_t)\bigr) - \bigl( f(\eta')-\EXPECT_{\eta_0} f(\eta_t) \bigr) \Bigr ]^2 \bar p(\sigma',\eta',t;\sigma_0,\eta_0)  \\
&= \sum_{\sigma'} \Bigl [ f(\sigma')-\EXPECT_{\sigma_0} f(\sigma_t) \Bigr ]^2 p^{A}(\sigma',t;\sigma_0) + \sum_{\eta'} \Bigl [f(\eta')-\EXPECT_{\eta_0}f(\eta_t) \bigr) \Bigr ]^2 p^{B}(\eta',t;\eta_0)  -  \\
&\quad\quad -2 \sum_{\sigma',\eta'}  \bigl(f(\sigma')-\EXPECT_{\sigma_0} f(\sigma_t)\bigr)\bigl( f(\eta')-\EXPECT_{\eta_0}f(\eta_t) \bigr) \bar p(\sigma',\eta',t;\sigma_0,\eta_0)    \PERIOD
\end{split}
\end{equation}
If these processes are independent then
\begin{equation}\label{joint_probability}
\bar p(\sigma',\eta',t;\sigma_0,\eta_0)  =  p^A(\sigma',t;\sigma_0)   p^{B}(\eta',t;\eta_0)
\end{equation}
 and the last term in equation \VIZ{variance_int} is 0. Hence a strong coupling of
the processes, i.e. when the processes $\sigma$ and $\eta$ are correlated, is expected to reduce the variance of the estimator (\ref{estimator_mean}). Therefore, based on (\ref{variance_int}), our goal is to maximize
\begin{equation}
\sum_{\sigma',\eta'} \bigl(f(\sigma')-\EXPECT_{\sigma_0} f(\sigma_t)\bigr)\bigl( f(\eta')-\EXPECT_{\eta_0}f(\eta_t) \bigr) \bar p(\sigma',\eta',t;\sigma_0,\eta_0)
\end{equation}
although note that inequality \VIZ{variance_int} implies that any coupling has a theoretical upper bound.


As shown above, highly correlated processes will reduce the variance of estimator \VIZ{estimator_mean}. Earlier work on coupling focuses on reactions networks which model well mixed systems, i.e. when the system is spatially homogeneous. The simplest
coupling is to run the two processes using the same stream of random numbers, known as Common Random Number (CRN). One improvement of this
method is the Common Random Path algorithm \cite{Khammash_CRP} (CRP). This
method uses  separate streams of random numbers for every reaction, introducing non zero correlation between the two processes thus leading
to more coupled processes. In systems where the
dynamics are high dimensional this coupling is intractable:  the number of reactions is proportional to the lattice
size which means that in order to implement the CRP method one has to store as many different random number streams or keep as many
different seeds and bring the random number generator to the previous state for every reaction. The former needs extremely large amount of
memory, while the later is using up too much computational time.

A different approach was introduced recently  where the two processes are coupled using a generator acting on both
processes \cite{Anderson}. The same idea was used by Liggett  \cite{Liggett}, in a different context, to obtain theoretical results for the monotonicity and ergodicity of interacting particle systems by Markov Chains. Our work 
uses these ideas as a starting point and extends them to spatial systems with complex dynamics and a high dimensional state space. 
More importantly our approach incorporates properties of macroscopic observables leading to a much more efficient goal-oriented method. Finally we also 
refer to Lindvall \cite{Lindvall} for the mathematical framework of coupling for random variables and general stochastic processes.

\section{Formulation of Coupling and Variance Reduction for Spatially Extended Systems}\label{section:formulation_of_coupling}

In this Section we first present couplings of stochastic processes based on correlating the  time series of the KMC evolutions defined on the entire configuration space. Subsequently, in Section 4 we present a class of more efficient couplings for variance reduction, which are nevertheless constructed building on the couplings of Section 3, and are based on correlating just the time series of the targeted observable. 

In the same spirit of earlier coupling methods, in order to reduce the variance of the estimator \VIZ{estimator_mean}, one can couple the two
generators and produce highly correlated paths. The novelty of our approach is that we introduce couplings that
are efficient for very high-dimensional state spaces arising in  spatially extended systems such as reaction-diffusion KMC.  We will present
a general and systematic formulation for the coupled generator and the rates of this generator will be obtained by solving an optimization
problem that is going to be presented subsequently, see Section \ref{sec:goal_oriented_optimization_problem} below.

We first define a coupling of stochastic process, $\zeta_t = (\PROCa,\PROCb)$, constructed so that (a) $\zeta_t$ is a continuous time Markov chain, hence it is easy to implement using KMC methods and (b) $\sigma_t$ (respectively $\eta_t$) is a Markov
process with generator $\LOPERa$ (respectively $\LOPERb$). We demonstrate the construction of such a process next. These  processes $\sigma_t,\eta_t$ need to
be strongly correlated as suggested by \VIZ{variance_int}. In order to achieve this we will define a generator for $\zeta_t$, see Section \ref{section:background_MC}, acting on an observable
function on both processes, that will couple the dynamics of both processes but at the same time it will respect the dynamics of each one.
More specifically, let $g(\sigma,\eta)$ be an observable on two processes. The corresponding coupled
generator
$\LOPERcpl$ should be related to generators $\LOPERa$ and $\LOPERb$ and must satisfy the following property
\begin{eqnarray}\label{cpl_property}
\LOPERcpl g(\sigma,\eta)=\LOPERa f(\sigma)\COMMA &\text{ if }& g(\sigma,\eta)=f(\sigma) \;\; \text{and} \nonumber  \\
\LOPERcpl g(\sigma,\eta)=\LOPERb f(\eta)   \COMMA &\text{ if }&   g(\sigma,\eta)=f(\eta)    \PERIOD
\end{eqnarray}
That is $\bar\LOPER$ reduces to $\LOPERa$ (respectively $\LOPERb$) for observables depending only on $\sigma$ (respectively $\eta$). Indeed, it
can be proved (see \cite{Liggett} Ch.3, Theorem 1.3) that if property \VIZ{cpl_property} holds then, 
\begin{eqnarray}\label{cpl_property_2}
\EXPECT_{(\sigma_0,\eta_0)}  g(\sigma_t,\eta_t) = \EXPECT_{\sigma_0 } f(\sigma_t)\COMMA &\text{ if }& g(\sigma,\eta)=f(\sigma) \;\; \text{and}
\nonumber  \\
\EXPECT_{(\sigma_0,\eta_0)}  g(\sigma_t,\eta_t) = \EXPECT_{\eta_0 } f(\eta_t)  \COMMA &\text{ if }&   g(\sigma,\eta)=f(\eta)    \COMMA
\end{eqnarray}
where $\mathbb{E}_{(\sigma_0,\eta_0)} g(\sigma_t,\eta_t)$ denotes the mean value of the observable function $g$ of the state variables
$(\sigma_t,\eta_t)$ with respect to the law imposed by the initial data $(\sigma_0,\eta_0)$. Heuristically this fact can also be seen, at least for short times, 
by combining (\ref{semigroup}), (\ref{approximation_semigroup}) and (\ref{cpl_property}). 
These relations imply  that all averages of two coupled processes, and thus all observables $f$, coincide with the averaged observables generated by the uncoupled generators $\LOPERa,\LOPERb$.


\subsection{Variance reduction via coupling}\label{sec:variance_reduction_via_coupling}
We are interested in the evolution of the variance of  an observable on two processes $(\sigma_t,\eta_t)$, imposed from the
dynamics of $\RATEa$ and $\RATEb$. This pair of processes can be simulated using $\LOPERa$ and $\LOPERb$ or the coupled generator
$\LOPERcpl$. Variance in (\ref{variance_int}) can be written as
\begin{equation}\label{variance}
\VAR\bigl [ f(\PROCa) -  f(\PROCb) \bigr ]  =\VAR f^2(\sigma_t) + \VAR f^2(\eta_t) + 2\EXPECT f(\sigma_t)\EXPECT f(\eta_t) - 2\EXPECT \bigl [ f(\sigma_t)f(\eta_t)\bigr] \COMMA
\end{equation}
where the mean and variance in the above equation is assumed with respect to initial data $(\sigma_0,\eta_0)$. From now on we will adopt
this lighter notation
unless otherwise is stated. The first three terms of the above equation cannot be controlled from coupling as they depend on the mean values
of single paths, which by \VIZ{cpl_property_2} are independent of coupling. However, the last term is the mean solution of a
coupled process $\zeta_t=(\sigma_t,\eta_t)$ with generator $\LOPERcpl$. Thus, if we pick $g(\sigma,\eta)=f(\sigma)f(\eta)$ and  define 
the mean value of $g$ with respect to initial data $(\sigma_0,\eta_0)$,
\begin{equation}\label{var_dep_obs}
u(t;\sigma_0,\eta_0)  =  \mathbb{E}_{(\sigma_0,\eta_0)}g(\PROCa,\PROCb) = \mathbb{E}_{(\sigma_0,\eta_0)}f(\PROCa)f(\PROCb)
\stackrel{\VIZ{semigroup}}{=} e^{t\LOPERcpl}f(\sigma_0)f(\eta_0)  \COMMA
\end{equation}
then $u$ satisfies the following differential equation,
\begin{equation}
\partial_t  u(t;\sigma_0,\eta_0) = \LOPERcpl u(t;\sigma_0,\eta_0), \quad   u(0;\sigma_0,\eta_0) = g(\sigma_0,\eta_0) \PERIOD
\end{equation}
The coupled generator $\bar \LOPER$ will be chosen in a way such that it maximizes the quantity in (\ref{var_dep_obs}) leading to a variance minimization in (\ref{variance}).

\subsection{Microscopic coupling}\label{sec:microscopic_coupling}

\begin{figure*}[htpb]
    \begin{center}
	  \includegraphics[scale=1.2]{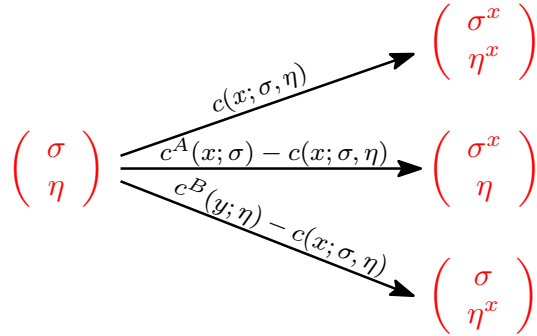}
	  \caption{The process being at state $(\sigma,\eta)$ can move to three different states for every $x\in\LATT$ according to the coupled generator (\ref{generator_cpl}). Each state has a different probability,
 $p_1(x;\sigma,\eta)=c(x;\sigma,\eta)/c_0(x;\sigma,\eta)$, $p_2(x,\sigma,\eta)=(c^A(x;\sigma)-c(x;\sigma,\eta))/c_0(x;\sigma,\eta)$, $p_3(y,\sigma,\eta)=(c^B(y;\eta)-c(x;\sigma,\eta))/c_0(x;\sigma,\eta)$ and $c_0 = c + c^A - c + c^B -c = c^A + c^B - c$.  }
\label{fig:transitions}
    \end{center}
\end{figure*}

For concreteness, we will first  present the idea of coupling using a simple Ising type model, see also (\ref{ising_rates}). We call this coupling a \textit{microscopic coupling} because it 
couples the time series of the entire configuration, i.e. it pairs every site $x$ of the $\sigma$ process with the same site $x$ of the $\eta$ process, see Figure \ref{fig:schematic_couples}. Further examples of typically more efficient couplings which are designed for specific observables will be presented in Section \ref{sec:coupling_macroscopic_observables}, see also Figure \ref{fig:schematic_couples}.  Returning to the introduction of microscopic couplings, the following coupling was introduced by Liggett  \cite{Liggett} to study monotonicity and ergodic properties of interacting particle systems,
\begin{equation}\label{generator_cpl}
\begin{split}
\LOPERcpl g(\sigma,\eta) =  &\sum_{x\in\LATT}  c(x;\sigma,\eta)  \Bigl [ g(\sigma^x,\eta^x) - g(\sigma,\eta) \Bigr ] \\
&+ \bigl ( \RATEa(x;\sigma) - c(x;\sigma,\eta) \bigr ) \Bigl [ g(\sigma^x,\eta) - g(\sigma,\eta) \Bigr ]  \\
&+ \bigl ( \RATEb(x;\eta)   - c(x;\sigma,\eta) \bigr ) \Bigl [ g(\sigma,\eta^x) - g(\sigma,\eta) \Bigr ] \COMMA
\end{split}
\end{equation}
for $\sigma(x),\eta(x) \in\Sigma=\{0,1\}$ where $\sigma^x$ is a flip of $\sigma(x)$ to $1-\sigma(x)$. Figure \ref{fig:transitions}
demonstrates this coupling by visualizing the possible transitions and demonstrates that the coupled process $\zeta_t$ can be simulated like any standard KMC mechanism. 
Note that the coupled generator (\ref{generator_cpl}) can be
written in terms of the original generators $\LOPERa$ and $\LOPERb$ as,
\begin{equation}\label{generator_cpl_2}
\begin{split}
\LOPERcpl g(\sigma,\eta) &= \LOPERa g(\sigma,\eta) + \LOPERb g(\sigma,\eta)  \\
 & \quad + \sum_{x\in\LATT} c(x;\sigma,\eta) \bigl [ g(\sigma^x,\eta^x) -g(\sigma^x,\eta) -g(\sigma,\eta^x) + g(\sigma,\eta) \bigr ] \PERIOD
\end{split}
\end{equation}
It is straightforward to verify that for all possible rates $c=c(x;\sigma,\eta)$ the above coupling satisfies property \VIZ{cpl_property}. Since $\zeta_t$ is Markovian, the rates in $\LOPERcpl$ (see also Figure \ref{fig:transitions}) must be non-negative, i.e., we have the condition
\begin{equation}\label{rate_property}
0  \leq  c(x;\sigma,\eta) \le  \min \{ \RATEa(x;\sigma) , \RATEb(x;\eta) \} \PERIOD
\end{equation}
The form of the rate function $c$ will be obtained next by solving an optimization problem that involves \VIZ{rate_property} as a
constraint while the variance is optimized.

In order to have a comparison basis, we define a new  generator, 
\begin{equation}\label{simplest_coupling}
\LOPERcplSimple g(\sigma,\eta) = \LOPERa g(\sigma,\eta) + \LOPERb g(\sigma,\eta)\COMMA
\end{equation} 
which corresponds to the case $c\equiv 0$ in (\ref{generator_cpl_2}). We call it the \textit{trivial coupling} because, as suggested by the construction in Figure \ref{fig:transitions}, there is no attempt to correlate the two time series $\sigma_t,\eta_t$, while \VIZ{cpl_property} is obviously satisfied. As indicated by the simulations in Section \ref{section:numerics_macro} this generator leads to same variance as with using two uncoupled processes, constructed by independent KMC.

We will compute the difference of the variance between the general coupled generator $\LOPERcpl$ that respects (\ref{cpl_property}) and the trivially coupled generator $\LOPERcplSimple$, showing that the variance of estimator \VIZ{estimator_mean} using the coupled generator is
always less or equal than that of using the trivially coupled generator. The computation is done locally in time, meaning that we examine the behavior of the variance for $t\in [0,\delta t ]$ for $\delta t \ll 1$.  Although the computations below follow these
constraints, the numerical results in Section \ref{section:numerics_micro} indicate that the theoretical result is also true for large times.  The presentation of a more general theoretical result is beyond the scope of this paper and can be considered by exploiting the mathematical tools presented in recent works \cite{AKP,KPS}.

In (\ref{variance}) we see that the variance of the coupled generator depends only on the quantity \VIZ{var_dep_obs}. Thus the
difference of the variance will depend on the same quantity,
\begin{equation}\label{semigroup_diff}
 \Bigl [ \EXP{\LOPERcpl \delta t} -  \EXP{\LOPERcplSimple \delta t} \Bigr ] g(\sigma_0,\eta_0)  = \delta t \bigl ( \LOPERcpl -
\LOPERcplSimple
\bigr) g(\sigma_0,\eta_0) + \BIGO (\delta t^2) \COMMA
\end{equation}
where we used \VIZ{solution_series} since both generators are bounded due to \VIZ{bounded}. The above expansion shows that the difference in
variance, locally in time, is controlled by the difference of the two generators. By comparing
 $\LOPERcplSimple$ with $\LOPERcpl$ from equation \VIZ{generator_cpl_2} and noting that $g(\sigma,\eta)=f(\sigma)f(\eta)$, the difference becomes,
\begin{equation}\label{oper_diff}
\begin{split}
&\bigl ( \LOPERcpl - \LOPERcplSimple \bigr)f(\sigma)f(\eta)  = \\
&=   \sum_{x\in\LATT} c(x;\sigma,\eta) \bigl [ f(\sigma^x)f(\eta^x)-f(\sigma^x)f(\eta) - f(\sigma)f(\eta^x) + f(\sigma)f(\eta) \bigr ] \\
&= \sum_{x\in\LATT} c(x;\sigma,\eta) \bigl [ f(\sigma^x) - f(\sigma) \bigr ] \bigl [ f(\eta^x) - f(\eta) \bigr ] \PERIOD
\end{split}
\end{equation}

\begin{remark}
If we set $\sigma_0=\eta_0$ in \VIZ{oper_diff}, which corresponds to the case where the two trajectories
have the same initial data, then,
\begin{equation}
\bigl ( \LOPERcpl - \LOPERcplSimple \bigr)f(\sigma_0)f(\sigma_0)  =   \sum_{x\in\LATT} c(x;\sigma_0,\sigma_0) \bigl (
f(\sigma_0^x)-f(\sigma_0) \bigr
)^2 \geq 0 \PERIOD
\end{equation}
This last relation, combined with (\ref{variance}),(\ref{var_dep_obs}) and (\ref{semigroup_diff}), clearly shows that the variance of $\LOPERcpl$ is always less or equal than that of the $\LOPERcplSimple$ for processes that start from the same state. 
\end{remark}

Next we show that we can obtain such variance reduction for processes starting at different initial configurations and determine the form of the coupled rates $c$ in (\ref{generator_cpl}) by defining and solving an appropriate optimization problem.

\subsection{Goal Oriented Optimization as a Coupling Design Tool}\label{sec:goal_oriented_optimization_problem} 

Indeed, in order to find the form and hopefully the optimal choice  of the rate function $c(x;\sigma,\eta)$ in (\ref{generator_cpl}), see also Figure \ref{fig:transitions},  we can write the variance reduction
problem as an optimization problem motivated by \VIZ{variance}. Note that minimization of variance is equivalent of maximizing (\ref{var_dep_obs}). Thus for $t=\delta t \ll 1$, equation (\ref{var_dep_obs})  can be rewritten as,
\begin{equation}
 \begin{split}
 \EXPECT \bigl [ f(\sigma_{\delta t})f(\eta_{\delta t})\bigr] & = \EXP{\LOPERcpl\delta t}f(\sigma_0)f(\eta_0) 
\nonumber \\
& =   f(\sigma_0)f(\eta_0) + \delta t \LOPERcpl f(\sigma_0)f(\eta_0)
+ \BIGO (\delta t^2)  \bigr ) \nonumber \PERIOD  
 \end{split}
\end{equation}
Using equation \VIZ{generator_cpl_2} we have:
\begin{equation}
\begin{split}
\EXPECT \bigl [ f(\sigma_{\delta t})f(\eta_{\delta t})\bigr] = & f(\sigma_0)f(\eta_0) + \delta t \bigl [ \LOPERa+\LOPERb \bigr
]f(\sigma_0)f(\eta_0)  \\ & +\delta t \sum_{x\in\LATT} c(x;\sigma_0,\sigma_0) \bigl [ f(\sigma_0^x)-f(\sigma_0) \bigr ]\bigl [
f(\eta_0^x)-f(\eta_0) \bigr ]+ \BIGO (\delta
t^2)
\PERIOD 
\end{split}
\end{equation}
Because the generator $\LOPERa$ acts only on the process $\sigma_t$ and $\LOPERb$ acts only on $\eta_t$ the above relation yields,
\begin{equation}\label{mean_product}
\begin{split}
\EXPECT \bigl [ f(\sigma_{\delta t})f(\eta_{\delta t})\bigr] = &f(\sigma_0)f(\eta_0) + \delta t f(\eta_0)\LOPERa f(\sigma_0) + \delta t
f(\sigma_0)\LOPERb f(\eta_0) \\ 
 & + \delta t \sum_{x\in\LATT} c(x;\sigma_0,\sigma_0) \bigl [ f(\sigma_0^x)-f(\sigma_0) \bigr ]\bigl [
f(\eta_0^x)-f(\eta_0) \bigr ]+ \BIGO (\delta t^2) \PERIOD
\end{split}
\end{equation}
As can be seen from the above equation, the only term we can control to reduce variance is the last term in last equation involving $c$. In order to
minimize the variance in one time step $\delta t$, we maximize the leading order of  (\ref{mean_product}) under the constraints \VIZ{rate_property}:

\begin{optProb}\label{optimization_problem_state_1}
Given the rate functions $c_A(x,\sigma)$ and $c_B(x,\eta)$ and  an observable function $f$, the rates $x(x;\sigma,\eta)$ of the coupled generator (\ref{generator_cpl}) that minimize the variance of estimator (\ref{estimator_mean}) are given by the solution of the optimization problem,
\begin{subequations}\label{max_problem}
\begin{align}
\max_{c}   \mathcal{F}[ c;f ] =  \max_{c} &  \sum_{x\in\LATT}   c(x;\sigma,\eta)  \bigl [ f(\sigma^x)-f(\sigma) \bigr ]\bigl [ f(\eta^x)-f(\eta ) \bigr ]  \label{max_problem_functional} \\
 &\text{under the constraint } \nonumber \\
  0  \leq  c(x;\sigma,&\eta)  \le  \min \{ \RATEa(x;\sigma) , \RATEb(x;\eta) \}  \label{max_problem_constraint}  \PERIOD
\end{align}
\end{subequations}
The constraint is needed to ensure positive rates in generator (\ref{generator_cpl}), see also Figure \ref{fig:transitions}.
\end{optProb}
The functional $\mathcal{F}$ will be used below as a diagnostic tool to design and evaluate different couplings for variance reduction.
One obvious, possibly suboptimal, choice for the coupled rates, satisfying the constraint in (\ref{max_problem_constraint}), is 
\begin{equation} \label{unoptimized_micro}
c_0(x;\sigma,\eta) = \min \{ \RATEa(x;\sigma) , \RATEb(x;\eta) \},
\end{equation}
which we will show next that it is not the optimal choice in this class of microscopic couplings. The maximization problem (\ref{max_problem}) depends clearly on the choice of the observable function $f$: if we choose, for example, as observable the coverage on the lattice, i.e.
\begin{equation}\label{coverage}
 f(\sigma) = \frac{1}{N} \sum_{x\in\LATT}\sigma(x) \COMMA
\end{equation}
then it is straightforward to verify that
\begin{equation}
 \SGN \Bigl ( \bigl [ f(\sigma^x)-f(\sigma) \bigr ]\bigl [ f(\eta^x)-f(\eta ) \bigr ]  \Bigr ) = 
 \begin{cases}
  +1, & \textrm{if } \sigma(x) = \eta(x) \\
  -1, & \textrm{if } \sigma(x) \neq \eta(x) 
  \end{cases} \PERIOD
\end{equation}
Thus the solution to the maximization problem \VIZ{max_problem} for observable (\ref{coverage}) is the function 
\begin{equation}\label{optimized_micro}
  c_1(x;\sigma,\eta) = 
  \begin{cases}
    \min \{ \RATEa(x;\sigma) , \RATEb(x;\eta) \},  & \textrm{if } \sigma(x) = \eta(x) \\
    0, & \textrm{if } \sigma(x) \neq \eta(x)
 \end{cases} \PERIOD
\end{equation}
which is exactly the generator proposed by Liggett \cite{Liggett}. Notice that both rates, (\ref{unoptimized_micro}) and (\ref{optimized_micro}), define microscopic couplings due to the fact that they couple the same site $x$ for the two processes through the mechanism depicted in Figure \ref{fig:transitions}.

\begin{figure*}[htpb]
	\begin{center}
	\includegraphics[scale=0.5]{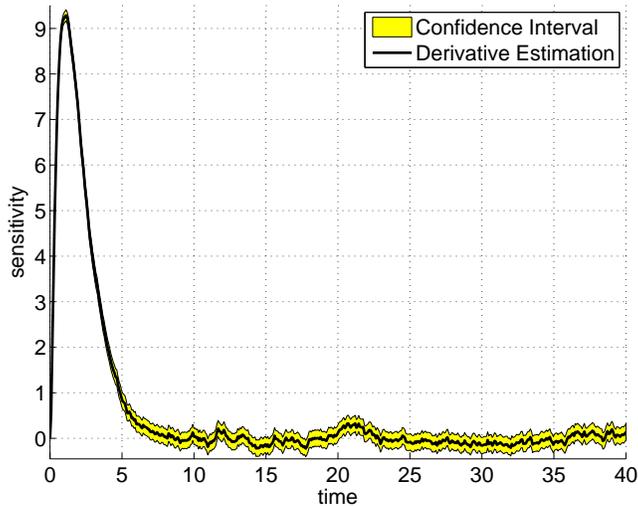}
	\caption{Estimation of derivative and confidence interval for the example discussed in Section \ref{section:numerics_micro} }	
	\label{fig:simulation1_derivative}	
	\end{center}
\end{figure*}

Finally for the functional (\ref{max_problem_functional}) holds that
\begin{equation}\label{first_inequality}
\mathcal{F}[c_0;f] \leq \mathcal{F}[c_1;f] \PERIOD
\end{equation}
Note also that from  (\ref{generator_cpl_2})  the rates $c\equiv 0$ is the case of the trivially coupled generator (\ref{simplest_coupling}).
\begin{equation}\label{second_inequality}
0 = \mathcal{F}[0;f] \leq \mathcal{F}[c_1;f] \PERIOD
\end{equation}
Therefore the ``tightest'' coupling is expected to be given by (\ref{optimized_micro}). Here, the functional $\mathcal{F}$ is used as a diagnostic tool to evaluate and design different couplings obtaining the one with the maximum variance reduction. We explore this issue next with a specific example.

\subsection{Numerical results}\label{section:numerics_micro}

\begin{figure*}[htpb]
	\begin{center}
	\includegraphics[scale=0.6]{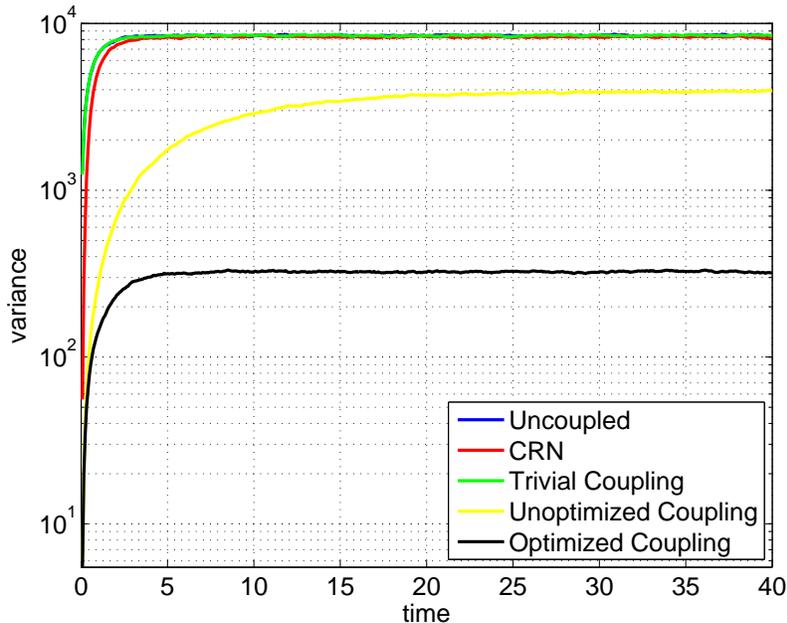}
	\caption{ Comparison of the variance for various estimation methods in Section \ref{sec:goal_oriented_optimization_problem}. The results are in agreement with inequality (\ref{first_inequality}) and (\ref{second_inequality}).}	
	\label{fig:simulation1_variance}	
	\end{center}
\end{figure*}

As a first example we will investigate the behavior of the nearest neighbor Ising model, with an adsorption/desorption mechanism described by the generator \VIZ{ising_generator}.
The states, $\sigma^x$, that the system can move to are given by (\ref{ising_new_configuration}) and the transition rate from $\sigma$ to $\sigma^x$ is given by (\ref{ising_rates}). 
For the simulation presented in Figure \ref{fig:simulation1_derivative} and \ref{fig:simulation1_variance} the parameters are $\beta=J=h=1$, the size of the lattice is $N=100$ and the final time
is $T=40$. In order to estimate the derivative with respect to $\beta$ we choose $\epsilon=0.1$ and the mean is computed over $4\times10^4$ sample paths.

In Figure \ref{fig:simulation1_derivative} the estimated derivative with respect to $\beta$, using the estimator (\ref{estimator_mean}) and coverage (\ref{coverage}) as observable,  is presented. The confidence interval of $99\%$ is also presented as an indicator that the calculated quantity is not far from the mean.
In Figure \ref{fig:simulation1_variance} the variance of various estimators for the simulation in Figure \ref{fig:simulation1_derivative} is presented. To be more precise, we compute the variance using  the uncoupled process, the CRN method, the trivial coupling (\ref{simplest_coupling}), the unoptimized coupling (\ref{unoptimized_micro}) and the optimized coupling (\ref{optimized_micro}). The first three schemes have variance of the same order which is about $10^4$. The unoptimized coupling  and the optimized coupling gives about an order and two orders variance reduction respectively. This calculation is in qualitative agreement with the theory presented in the previous section and specifically with inequality (\ref{inequality_1N}). Moreover, the results presented in this Section are all reproducible by Matlab scripts\cite{matlab_code}. For a short description of the implementation  see Appendix \ref{app:implementation}.

Note that in this simulation the sampling of the finite difference estimator (\ref{estimator_mean})  using the optimal coupling (\ref{optimized_micro}) requires $200$ times  fewer samples than the CRN or the uncoupled methods. At the same time the computational overhead of using the specific coupling slows down the KMC only by a factor of two leading to a scheme that is $100$ times faster than schemes used up to now such as CRN.


\section{Coupling Macroscopic Observables}\label{sec:coupling_macroscopic_observables}

Microscopic couplings defined in Section \ref{section:formulation_of_coupling}, at least heuristically, appear to be  too restrictive because the same site of the first lattice is coupled with the same site of the second lattice, see Figure \ref{fig:coupling_1_vs_N}. It is reasonable to conjecture that a  more flexible approach could be to couple certain observables of the system: for example, if we are interested in the concentration of certain species it would be preferable to couple the two systems in such a way that they will lead to correlated  concentrations rather than trying to correlate the {\em entire} configurations of the two systems. We call this class of couplings \textit{macroscopic couplings} due to the fact that they depend on macroscopic quantities that depend on the entire  lattice, e.g. coverage, Hamiltonian, spatial correlations, surface roughness, etc, in contrast to microscopic couplings that depend on local spin configurations.

First, we will present the formulation of  the generator corresponding to  a general lattice KMC model, extending the discussion of the previous Section; subsequently we will present the concept of macroscopic  coupling. 
A general optimization problem will be defined in Section  \ref{sec:macroscopic_coupling} giving the ability to compare this new  more general coupling with the microscopic coupling already defined
in Section~\ref{sec:microscopic_coupling}. Finally, we will embed all these couplings within a hierarchy of couplings (see Figure \ref{fig:schematic_couples}) defined on a hierarchy of {\em mesoscopic} geometric decompositions of the lattice, see also Figure \ref{fig:schematic_couples}.

\subsection{General form of the Generator}\label{section:general_form_of_the_generator}

\begin{figure*}[htpb]
\centering
	\includegraphics[scale=0.75]{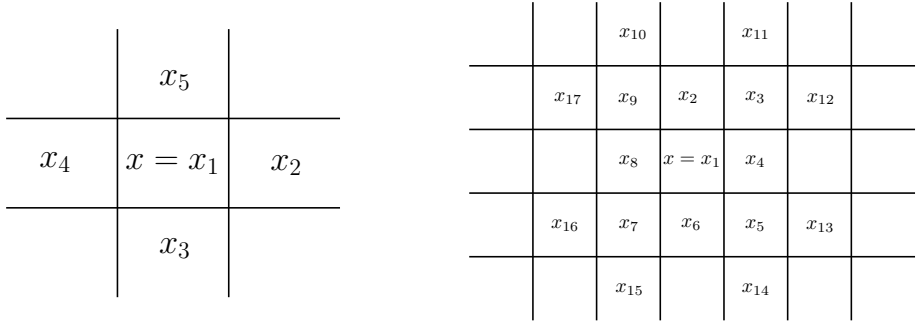}
      \caption{  Numbering of the neighborhood $\NEIGH_x$  for (left) the Ising or the ZGB \cite{ZGB} model and (right) a larger neighborhood of a more complex model of catalytic CO oxidation (see \cite{Evans}) (right), see also Appendix \ref{app:examples}.}
\label{fig:neighborhood}
\end{figure*}

In order to describe a lattice KMC model in great generality one needs to define and enumerate the set of all involved species, $\Sigma$, and how these species
interact with each other. The interactions are defined through the rate functions, $c(\sigma,\sigma')$, see (\ref{totalrate}),(\ref{skeleton}),(\ref{generator}). These functions give the
rate at which the system at given state $\sigma$, will move to a different state $\sigma'$ according to (\ref{totalrate}),(\ref{skeleton}). These changes are generally local, meaning that
affect only the sites in the neighborhood of a given site $x\in\LATT$, and thus the rate functions are also local functions. In order to make clear this idea we will introduce some notation and then present the generator specifics.

\begin{definition}\label{def:neigh}
 For a given $x\in\LATT$ define $\NEIGH_x =( x_1=x,x_2, \ldots, x_k )$, an ordering of \textbf{the neighborhood of site $x$} which
describes the local geometry, see Figure \ref{fig:neighborhood}. The state variable evaluated at the neighborhood sites will be denoted by
$\sigma_i=\sigma(x_i),\;i=1,\ldots,k$ The neighborhood is fully determined by the dynamics of the selected model. 
\end{definition}

\begin{definition}\label{def:local_conf_omega}
 Given a neighborhood $\NEIGH_x$ of size $k$ we define the \textbf{local configuration} as a vector
$\omega=(\omega_0,\ldots,\omega_k)\in\Sigma^k$, where $\Sigma^k$ is the set of $k$-tuples with elements in $\Sigma$.
The \textbf{global updated configuration} { affected by the local configuration $\omega$ at site $x$} is
defined as,
\begin{equation}
\CONFNEW(y) = 
\begin{cases}
 \omega_i, &    y=x_i\in\NEIGH_x,\; i=1,\ldots,k \\
 \sigma(y), & y \notin \NEIGH_x
\end{cases}  \PERIOD
\end{equation}
and the rate going from $\sigma$ to $\sigma^{x,\omega}$ is $c(x;\sigma,\omega)$.
Moreover we define the \textbf{set of all accessible local configuration} for the given dynamics of the system, i.e., the set of all local configurations such that for a given configuration $\sigma$ there exists at least one $x\in\LATT$ with non zero rate, 
\begin{equation}
\CONF(\sigma) := \{ \omega\in\Sigma^k \;|\; \textrm{ exists } x \in \LATT  \textrm{ such that } c(x,\omega;\sigma)\neq 0 \} \PERIOD
\end{equation}
The dependence of $\Omega$ on the configuration $\sigma$ will be omitted in the sake of a lighter notation. 
\end{definition}
Using these definitions the generator can be written in the following way,
\begin{equation}\label{generator_x_omega}
  \LOPER f(\sigma) = \sum_{x\in\LATT}\sum_{\omega\in\CONF } c(x,\omega;\sigma) \bigl [ f(\CONFNEW) - f(\sigma) \bigr ] \COMMA
\end{equation}
which will allow us to give the formulation for a general coupling. All spatial KMC models that have mechanisms such as adsorption, desorption, diffusion and reaction between species, 
see Section \ref{sec:lattice_kmc} and Appendix \ref{app:examples}, and much more complex KMC models \cite{SVgraph}, where molecules have internal degrees of freedom and may occupy a neighborhood on a lattice,  can be described by the generator (\ref{generator_x_omega}), or combinations of such generators.

The starting point before we describe the macroscopic coupling will be to present the microscopic coupled generator in its full generality as it corresponds to (\ref{generator_x_omega}). In analogy to (\ref{generator_cpl}) we couple every pair  $(x,\omega)$ of the first process, $\sigma_t$, with every pair $(y,\omega')$ of the second process, $\eta_t$. The following generator describes this concept,
\begin{equation}
\begin{split} \label{general_coupled}
\LOPERcpl & g(\sigma,\eta) = \sum_{x\in\LATT} \sum_{y\in\LATT}\sum_{\omega\in \CONF}\sum_{\omega' \in \CONF} c(x,y,\omega,\omega';\sigma,\eta) \bigl [   g(\sigma^{x,\omega},\eta^{y,\omega'})-g(\sigma,\eta) \bigr ] \\
&+ \sum_{x\in\LATT} \sum_{\omega\in \CONF} \Bigl [ c_A(x,\omega;\sigma) - \sum_{y\in\LATT} \sum_{\omega' \in \CONF}  c(x,y,\omega,\omega';\sigma,\eta) \Bigr ] \bigl [   g(\sigma^{x,\omega},\eta)-g(\sigma,\eta) \bigr ] \\
&+ \sum_{y\in\LATT}\sum_{\omega' \in \CONF} \Bigl [ c_B(y,\omega';\eta)  - \sum_{x\in\LATT} \sum_{\omega\in \CONF} c(x,y,\omega,\omega';\sigma,\eta) \Bigr ] \bigl [   g(\sigma^{x,\omega},\eta^{y,\omega'})-g(\sigma,\eta) \bigr ] \COMMA
\end{split}
\end{equation}
where the coupling rates $c(x,y,\omega,\omega';\sigma,\eta)$ should be defined in  a way such  that all rates are positive (see the constraints in \VIZ{max_problem_general}).
The microscopic coupling (\ref{generator_x_omega}) is a generalization of the microscopic coupling in Figure \ref{fig:transitions}.

\subsection{Macroscopic coupling}\label{sec:macroscopic_coupling}

In this section we will describe a different approach for coupling two processes $(\sigma_t,\eta_t)$ that is based on partitioning the range of an observable function. Instead of coupling the events that happen on each lattice site, as in the microscopic coupling (\ref{generator_cpl}), we choose to couple events that belong to the same predefined class determined by level sets of the targeted observable $f$. This is natural choice since we are interested in the sensitivity analysis of a specific observable, hence we need to couple only the time series of that observable in (\ref{estimator_mean}). In turn, we expect that such a perspective can allow us more flexibility in the choice of coupling methods in order to minimize the variance of the sensitivity estimators. Indeed, as we demonstrate below, this goal-oriented approach will lead to greater variance reduction compared to the microscopic coupling. 


First, we define a partition of the range of the discrete derivatives, $f(\sigma^{x,\omega})-f(\sigma)$, into a finite number of classes, for example all level sets of $f(\sigma^{x,\omega})-f(\sigma)=k$.  The reason is that the macroscopic coupling will be based on coupling pairs of $(x,\omega)$ with $(y,\omega')$ that belong to the same class. 
We are interested in correlating time series of the observable $f$, hence we need to keep track of changes in $f$.
This partition depends on the choice of the observable function $f$.

\begin{figure*}[htpb]
     \centering
     \includegraphics[scale =0.6 ]{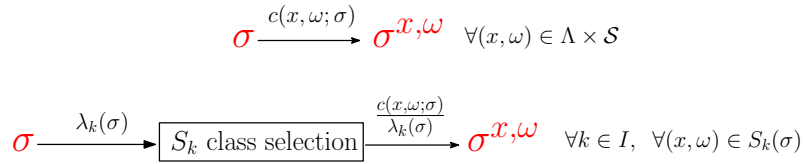}
	 \caption{ Two different implementations of KMC for the generator (\ref{generator_x_omega}) or equivalently (\ref{generator_rewriten}). In the top case the choice is being done in a single step with the appropriate rate, see (\ref{generator_x_omega}). In the bottom case the transition is done in two steps: first choose one of the predefined sets $S_k$ and then choose a state from that set using a new normalized rate, see (\ref{generator_rewriten}).}
	 \label{fig:ssa}
\end{figure*}

\begin{definition}\label{def:S_k_sets}
We decompose the range of an observable function $f$ in disjoint sets, i.e., let $J_i \subset \mathbb{R}, \;i\in I$ with $J_i \cap J_j = \emptyset$ when $i\neq j$ and $\cup_{i\in I}J_i = \mathbb{R}$. We now define the sets $S_k(\sigma)$ containing all possible events $(x,\omega)$ such that the value of $f(\sigma^{x,\omega})-f(\sigma)$ belongs in $J_k$, 
\begin{equation}
S_k(\sigma) = S_k(\sigma;f) = \bigl \{  (x,\omega) \in \LATT \times \CONF \;\;|\;\;  c(x,\omega;\sigma)\neq 0 \;\; \textrm{and} \;\;  f(\sigma^{x,\omega})-f(\sigma) \in J_k  \bigr \} \PERIOD 
\end{equation}
\end{definition}
For example, one choice of partition are the sets
\begin{equation}\label{coverage_sets}
J_1=(-\infty,0),\quad  J_2=\{0\}, \quad J_3 = (0,\infty) \PERIOD
\end{equation} 
which separate the values of $f(\sigma^{x,\omega})-f(\sigma)$ according to their sign only.
We also define the corresponding rate, 
\begin{equation}\label{sk_rate}
\lambda_k(\sigma) = \sum_{(x,\omega)\in S_k(\sigma)} c(x,\omega;\sigma) \COMMA
\end{equation}
which is the total rate in which the set $S_k(\sigma)$ can be selected, see Figure \ref{fig:ssa}. Thus the probability of selecting a pair  $(x,\omega)$ from $S_k(\sigma)$ is ${c(x,\omega;\sigma)}/{\lambda_k(\sigma)}$.
Using these definitions we can rewrite the generator for a single process  (\ref{generator_x_omega}), as
\begin{equation}\label{generator_rewriten}
\begin{split}
\LOPER  f (\sigma) &= \sum_{x\in\LATT} \sum_{\omega\in \CONF } c(x,\omega;\sigma)\bigl [   f(\sigma^{x,\omega})-f(\sigma) \bigr ] \\
&= \sum_{k\in I} \sum_{(x,\omega) \in S_k(\sigma)} c(x,\omega;\sigma)\bigl [   f(\sigma^{x,\omega})-f(\sigma) \bigr ] \\
&= \sum_{k\in I}  \lambda_k(\sigma)  \sum_{(x,\omega) \in S_k(\sigma)}  \frac{c(x,\omega;\sigma)}{\lambda_k(\sigma)}   \bigl [   f(\sigma^{x,\omega})-f(\sigma) \bigr ] \PERIOD
\end{split}
\end{equation}
Note that  in the straightforward implementation (top of Figure \ref{fig:ssa}, corresponding to first line of (\ref{generator_rewriten})) only one transition is needed. In our case (bottom Figure \ref{fig:ssa}, bottom line of (\ref{generator_rewriten})) we first choose one of the predefined sets $S_k(\sigma)$ and then choose a pair $(x,\omega)$ from this set. Obviously these two methods are equivalent because they give rise to the same stochastic process since they have the same generator \cite{Liggett}.

\begin{figure*}[htpb]
     \centering
     \includegraphics[scale =0.6 ]{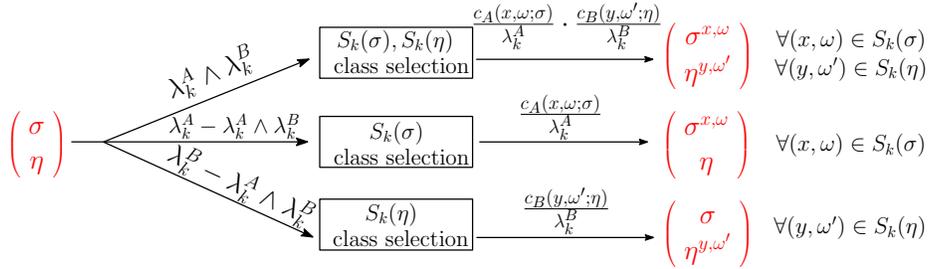}
	 \caption{ Schematic representation of the macroscopically coupled generator (\ref{macro_coupled}).  The coupling involves only the selection of the set $S_k(\sigma)$, see Definition \ref{def:S_k_sets}, and/or $S_k(\eta)$ (first transition) and then a state from each one of these  sets is selected independently (second transition). Compare to microscopic coupling in Figure \ref{fig:transitions}.}
	 \label{fig:ssa_coupled}
\end{figure*}

\begin{remark}\label{remark:BKL}
The Bortz-Kalos-Lebowitz algorithm \cite{BKL}, also known as the 10-fold algorithm, can be written in the form of equation (\ref{generator_rewriten}), see also Figure \ref{fig:ssa_coupled}. In the BKL algorithm the lattice sites of the Ising model (\ref{ising_rates}) are grouped in ten sets of equal rate. Thus, the ``classes'' $S_k(\sigma)$ are defined through
\begin{equation*}
    S_k(\sigma) = \{ x\in\LATT \; | \; c(x;\sigma) = c_k \}, \quad k = 1,\ldots,10  \COMMA
\end{equation*}
where $c_k,\;k=1,\ldots,5$ is the rate at a site with spin $0$ having $0$ to $4$ neighbours and  $k=6,\ldots,10$ is the rate at a site with spin $1$ having $0$ to $4$ neighbours. The generator of the Ising process is written as,
\begin{equation*}
\LOPER f(\sigma) = \sum_{k} c_k  \sum_{x\in S_k(\sigma)} \frac{c(x;\sigma)}{c_k}   \bigl [   f(\sigma^{x})-f(\sigma) \bigr ] \PERIOD
\end{equation*}
\end{remark}

\noindent \textbf{Macroscopic Coupling.}   Equipped with the equivalent representation (\ref{generator_rewriten}), we turn our attention to macroscopic coupling.
In analogy with definition (\ref{sk_rate}) we define the rates for the two processes $\sigma$ and $\eta$,
\begin{equation}
\lambda_k^A(\sigma) = \sum_{(x,\omega)\in S_k(\sigma)} c_A(x,\omega;\sigma) \quad \textrm{and} \quad \lambda_k^B(\eta) = \sum_{(y,\omega')\in S_k(\eta)}  c_B(y,\omega';\eta) \PERIOD
\end{equation}
Instead of coupling every lattice site $x$ 
as in Figure \ref{fig:transitions}, in this approach the coupling is being done in the selection of $S_k(\sigma)$ and $S_k(\eta)$ and then a state from each of these sets is selected independently
from each other. Thus the rate at which the set $S_k(\sigma)$ and/or $S_k(\eta)$ is being chosen has three steps:
\begin{enumerate}
 \item with the minimum of $\lambda_k^A(\sigma)$ and $\lambda_k^B(\eta)$ we choose the sets $S_k(\sigma)$ and $S_k(\eta)$, i.e., do a transition in both processes, $\sigma$ and $\eta$, 
 \item with $\lambda_k^A(\sigma)$ subtracted the minimum from step 1. we choose the set $S_k(\sigma)$, i.e., move only the $\sigma$ process  and
 \item with $\lambda_k^B(\eta)$      subtracted the minimum from step 1. we choose the set $S_k(\eta)$     , i.e., move only the $\eta$ process.
\end{enumerate}
According to this mechanism, after selecting the sets $S_k(\sigma)$ and/or $S_k(\eta)$, a pair $(x,\omega)$ from $S_k(\sigma)$ is selected independently from  a pair $(y,\omega')$ from $S_k(\eta)$.

In Figure \ref{fig:ssa_coupled} a schematic representation of this procedure is presented. The generator that describes this coupling in the spirit of (\ref{generator_rewriten}) is:
\begin{equation}
\begin{split} \label{macro_coupled}
&\LOPERcpl g(\sigma,\eta)  = \sum_{k\in I}  \min\{ \lambda_k^A(\sigma) , \lambda_k^B(\eta) \} \;\; \ldots \\
& \qquad\qquad \ldots\sum_{(x,\omega) \in S_k(\sigma)} \frac{c_A(x,\omega;\sigma)}{\lambda_k^A(\sigma)}\sum_{(y,\omega') \in S_k(\eta)} \frac{c_B(y,\omega';\eta)}{\lambda_k^B(\eta)} \bigl [   g(\sigma^{x,\omega},\eta^{y,\omega'})-g(\sigma,\eta) \bigr ] \\
&+ \sum_{k\in I} \Bigl [   \lambda_k^A(\sigma) -  \min\{ \lambda_k^A(\sigma) , \lambda_k^B(\eta) \}  \Bigr]  \sum_{(x,\omega) \in S_k(\sigma)} \frac{c_A(x,\omega;\sigma)}{\lambda_k^A(\sigma)} \bigl [   g(\sigma^{x,\omega},\eta)-g(\sigma,\eta) \bigr ] \\
&+ \sum_{k\in I} \Bigl [   \lambda_k^B(\eta) -  \min\{ \lambda_k^A(\sigma) , \lambda_k^B(\eta) \}  \Bigr] \sum_{(y,\omega') \in S_k(\eta)} \frac{c_B(y,\omega';\eta)}{\lambda_k^B(\eta)} \bigl [   g(\sigma,\eta^{y,\omega'})-g(\sigma,\eta) \bigr ] \PERIOD
\end{split} 
\end{equation}
A straightforward calculation shows that property \VIZ{cpl_property} of a coupling generator is satisfied. Moreover, all rates in (\ref{macro_coupled}), depicted also in Figure \ref{fig:ssa_coupled}, are positive functions.
 Furthermore, reordering the terms in \VIZ{macro_coupled} we can also write the total coupled rate for the general coupled generator \VIZ{general_coupled},
\begin{equation}\label{macro_coupled_rates}
\begin{split}
&c_N(x,y,\omega,\omega';\sigma,\eta) = \\
& \qquad\qquad =\sum_{k\in I}  \min\{ \lambda_k^A(\sigma) , \lambda_k^B(\eta) \}   \frac{c_A(x,\omega;\sigma)}{\lambda_k^A(\sigma)}  \frac{c_B(y,\omega';\eta)}{\lambda_k^B(\eta)} \chi_{S_k(\sigma)}(x,\omega) \chi_{S_k(\eta)}(y,\omega')  \COMMA
\end{split}
\end{equation}
where
\begin{equation}\label{characteristic}
\chi_A(a) = 
	\begin{cases}
			1,\quad \textrm{if } a\in A \\
			0,\quad \textrm{otherwise} \COMMA
	\end{cases}
\end{equation}
is the \textit{characteristic function} of the set $A$.
We will use this expression later in order to compare the gain in variance reduction of the macroscopic against the microscopic coupling.

\begin{figure*}[htbp]
    \begin{center}
	  \includegraphics[scale=0.4]{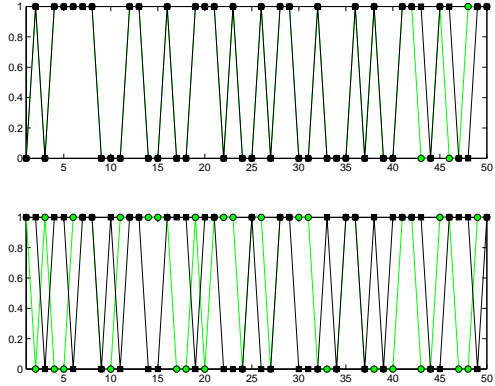}
	  \caption{Microscopic vs macroscopic snapshots at time $T=100$ from the coupled  Ising model using the generator (\ref{generator_cpl}) (top) and the generator (\ref{macro_coupled}) (bottom). The black circles and the green squares represent the two different systems for different parameters. Note that although in the first case there is very good pointwise agreement between the two   
	  configurations  while for the second this is not true, the observable used for this simulation (coverage) has values that are very close for both couplings. However, the corresponding variances are different, see Figure \ref{fig:var_vs_q_vs_t}. }
\label{fig:coupling_1_vs_N}
    \end{center}
\end{figure*}

\begin{figure*}[t]
                \begin{center}
                \includegraphics[scale=0.6]{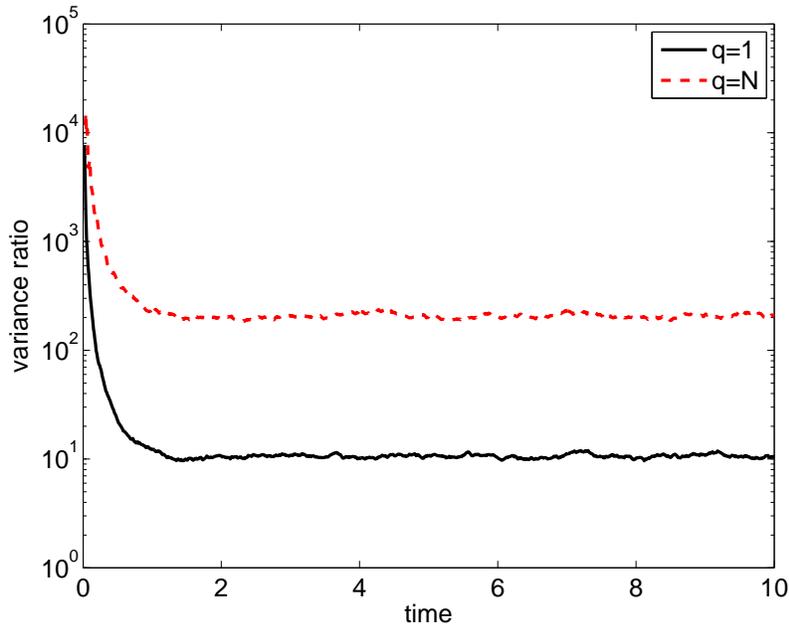}
        \caption{   Comparison of the variance ratio  between the microscopic coupling (\ref{generator_cpl}) $(q=1)$ and the macroscopic coupling (\ref{macro_coupled})$(q=N)$ for the 								 adsorption/desorption/diffusion model discussed in Section \ref{section:numerics_macro}. }
		\label{fig:var_vs_q_vs_t}     	
     	\end{center}
\end{figure*}

In Figure \ref{fig:coupling_1_vs_N} a pathwise comparison of the microscopic coupling (\ref{generator_cpl}) versus the macroscopic coupling (\ref{macro_coupled}) is presented for the example discussed in Section \ref{section:numerics_macro}. Note that in the first snapshot there is almost perfect pointwise agreement between the two simulated configurations
 while in the second snapshot the distributions of particles is different between the two systems in both couplings depicted in Figure \ref{fig:coupling_1_vs_N}. However, the value of the observable function, which is coverage in this example, is  0.5 for the $\sigma$ and 0.52 for the $\eta$ process in both depicted couplings in Figure \ref{fig:coupling_1_vs_N}. 
This fact will be discussed in more details in the next section where we show that the values of the observable function are tighter coupled for the macroscopic coupling (\ref{macro_coupled_rates}), leading to a greater variance reduction, see Figure \ref{fig:var_vs_q_vs_t}.

\subsection{Comparison of Microscopic and Macroscopic Couplings}\label{section:general_optimization_problem}

In the same spirit as in Section \ref{sec:goal_oriented_optimization_problem} we can formulate the maximization problem (\ref{max_problem}) for the general coupled generator \VIZ{general_coupled}. 
The functional $\mathcal{F}$ is obtained by the same procedure as the one presented in Section \ref{sec:goal_oriented_optimization_problem}. The analogue of problem \VIZ{max_problem} is stated bellow.

\begin{optProb}\label{optimization_problem_state_2}
Given a collection of rate functions $c_A(x;\omega,\sigma)$ and $c_B(x;\omega,\eta)$ and  an observable function $f$, we define
\begin{equation}\label{optimization_functional_2}
\mathcal{F} [c; f] = \sum_{x\in\LATT} \sum_{y\in\LATT}  \sum_{\omega\in \CONF}\sum_{\omega' \in \CONF}  c(x, y, \omega, \omega', \sigma,\eta)   [f(\sigma^{x,\omega})-f(\sigma)][f(\sigma^{y, \omega'})-f(\sigma)] \PERIOD
\end{equation}
Then  the rates $c(x, y, \omega, \omega', \sigma,\eta)$ of the coupled generator (\ref{general_coupled}) that minimize the variance of estimator (\ref{estimator_mean}) are given by the solution of the optimization problem,
\begin{subequations}\label{max_problem_general}
\begin{align}
&\max_c \mathcal{F} [c; f] \COMMA   \\
 \text{under the constraints } & \;\; 0  \leq  \sum_{y\in\LATT} \sum_{\omega' \in S}  c(x,y,\omega,\omega';\sigma,\eta)   \le  c_A(x,\omega;\sigma)  \COMMA \\
                                              & \;\; 0  \leq  \sum_{x\in\LATT} \sum_{\omega \in S}  c(x,y,\omega,\omega';\sigma,\eta)   \le  c_B(y,\omega';\eta)  \PERIOD
\end{align}
\end{subequations}
The constraints stem from the requirement that rates in generator \VIZ{general_coupled} are positive. 
\end{optProb}

It is straightforward to verify that rates \VIZ{macro_coupled_rates} of the macroscopically coupled generator \VIZ{macro_coupled} satisfy the constraints of optimization problem \VIZ{max_problem}. Thus, this coupling is a possibly suboptimal solution of the optimization problem (\ref{max_problem_general}). 

\begin{figure*}[htpb]
		\begin{center}	
	    \includegraphics[scale=0.6]{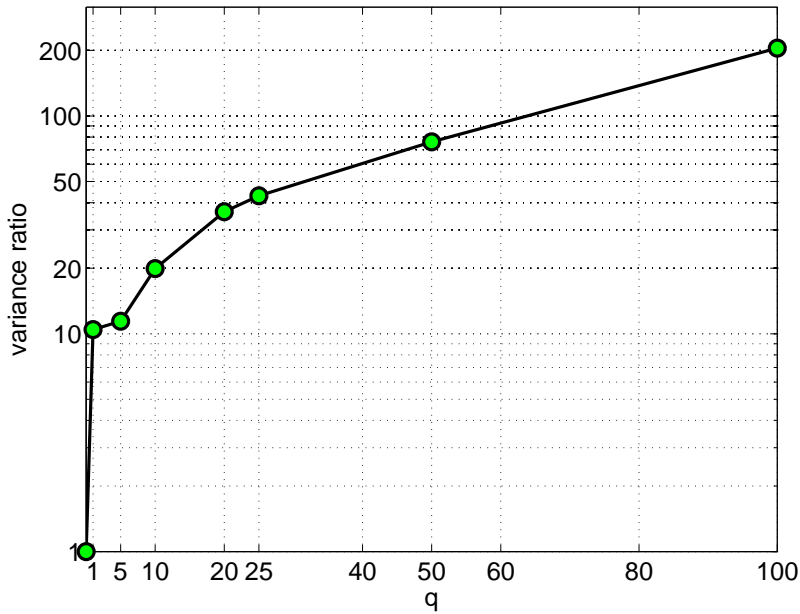}
        \caption{  Numerical verification of inequality (\ref{inequality_1qN}) where coverage (\ref{coverage}) is used as the observable. The first point in the graph $(q=0)$ corresponds to             						the uncoupled generator (\ref{simplest_coupling}). Note that the optimal choice of $q$ is obtained for $q=N$, at least  for the observable (\ref{coverage}). }
	 \label{fig:var_vs_q}
	 \end{center}
\end{figure*}

Our next goal is to compare the macroscopic coupling (\ref{macro_coupled_rates}) with the microscopic couplings (\ref{unoptimized_micro}),(\ref{optimized_micro}) presented in Section \ref{sec:goal_oriented_optimization_problem}.
We use the maximization functional (\ref{optimization_functional_2}) in order to compare different couplings, following along the lines of the simpler setting in section \ref{sec:goal_oriented_optimization_problem}. 
In fact, it is easy to see that  these microscopically coupled generators can be written in the general form of equation \VIZ{general_coupled}. For example, the rates for the unoptimized microscopic coupling \VIZ{unoptimized_micro}  using coverage as an observable are,
\begin{equation}\label{micro_rates_1}
 c_0(x,y,\omega,\omega';\sigma,\eta) =\min \{ c_A(x,\omega;\sigma) , c_B(y,\omega';\eta) \} \delta(x-y)  \COMMA
\end{equation}
and for the optimized microscopic coupling \VIZ{optimized_micro} ,
\begin{equation}\label{micro_rate_2}
\begin{split}
& c_1(x,y,\omega,\omega';\sigma,\eta) = \\
& \qquad\qquad = \min \{ c_A(x,\omega;\sigma) , c_B(y,\omega';\eta) \} \delta(x-y) \sum_{k\in I} \chi_{S_k(\sigma)}(x,\omega)\chi_{S_k(\eta)}(y,\omega') \COMMA
\end{split}
\end{equation}
where $J_k,\; k=1,2,3$ are the sets defined in (\ref{coverage_sets}) and $\chi$ the characteristic function (\ref{characteristic}).  As we already know from (\ref{first_inequality}),
\begin{equation}
\mathcal{F}[c_0;f] \leq \mathcal{F}[c_1;f] \COMMA
\end{equation}
where $f$ has been chosen to be the coverage defined by  \VIZ{coverage}, 
$$f(\sigma)=\frac{1}{N} \sum_{x\in\LATT} \sigma(x) \PERIOD$$
If we compare the expressions of $c_1$ and $c_N$ (\ref{macro_coupled_rates}) we can deduce that 
\begin{equation}\label{inequality_1N}
\mathcal{F}[c_1;f] \leq \mathcal{F}[c_N;f] \COMMA
\end{equation}
showing that the macroscopic coupling \VIZ{macro_coupled_rates} may not be an optimal solution since we did not optimize the functional (\ref{optimization_functional_2}) but reduces the variance compared to the previous microscopic couplings (\ref{unoptimized_micro}),(\ref{optimized_micro}). This inequality is a special case of inequality (\ref{inequality_1qN}) which will be proved in Appendix \ref{appendix:inequality_proof}. In Section \ref{section:numerics_macro} a numerical demonstration  of this inequality is presented, see Figure \ref{fig:var_vs_q}.
\begin{figure*}[htpb]
\centering
    	\includegraphics[scale=0.75]{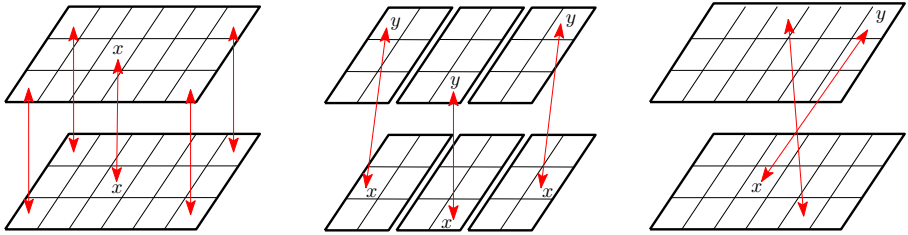}
      \caption{Schematic representation of the couplings,  microscopic (left), coarse grained (middle) and macroscopic (right).}
 	  \label{fig:schematic_couples}
\end{figure*}

\noindent \textbf{ Hierarchy of Couplings.} It is natural to ask the question if there are any other couplings between the finest microscopic and the coarsest macroscopic level. Indeed, such a coupling is a natural construction coming from coarse graining ideas \cite{KMV}: the lattice is split into $M$ cells, denoted $C_i,\;i=1,\ldots,M$ and the size of every cell is $|C_i|=q$. The idea is that instead of coupling the configurations $(x,\omega),(y,\omega')$ for $x,y$ in the whole lattice, as in the macroscopic coupling, we couple only the configurations for $x,y$ that belong in the same cell indexed by $i$ bellow, see Figure \ref{fig:schematic_couples}. The coupled rates in this case take the form,
\begin{equation}\label{coarse_coupled_rates}
\begin{split}
&c_q(x,y,\omega,\omega';\sigma,\eta) = \\ 
&=\sum_{i=1}^M  \sum_{k\in I}  \min\{ \lambda_{k,i}^A(\sigma) , \lambda_{k,i}^B(\eta) \}   \frac{c_A(x,\omega;\sigma)}{\lambda_{k,i}^A(\sigma)}  \frac{c_B(y,\omega';\eta)}{\lambda_{k,i}^B(\eta)} \chi_{ S_{k,i}(\sigma) } (x,\omega) \chi_{S_{k,i}(\eta)}(y,\omega')  \COMMA
\end{split}
\end{equation}
where
\begin{equation}
S_{k,i}(\sigma;f) = \bigl \{  (x,\omega) \in C_i \times \CONF \;\;|\;\;  c(x,\omega;\sigma)\neq 0 \;\; \textrm{and} \;\;  f(\sigma^{x,\omega})-f(\sigma) \in J_k  \bigr \} \COMMA 
\end{equation}
and
\begin{equation}\
\lambda_{k,i}^A(\sigma) = \sum_{(x,\omega)\in S_{k,i}(\sigma)} c_A(x,\omega;\sigma) \quad \textrm{and} \quad \lambda_{k,i}^B(\eta) = \sum_{(y,\omega')\in S_{k,i}(\eta)}  c_B(y,\omega';\eta) \PERIOD
\end{equation}
with the convention that $S_{k,N} = S_k$. Note that when $M=1$, i.e. $q=N$, then the above rates are the same as the rates in \VIZ{macro_coupled_rates} and when $M=N$, i.e., $q=1$ then the above rates are equal to the rates of equation (\ref{micro_rate_2}). We can use the functional (\ref{optimization_functional_2}) to compare these types of couplings. Indeed, for $f$ defined by (\ref{coverage}) we can show that (see Appendix \ref{appendix:inequality_proof}),
\begin{equation}\label{inequality_1qN}
\mathcal{F}[c_0;f]   \leq \mathcal{F}[c_1;f]   \leq \mathcal{F}[c_q;f]    \leq \mathcal{F}[c_N;f], \quad 1<q<N   \PERIOD
\end{equation}
Note that, as in inequality (\ref{inequality_1N}), the functional $\mathcal{F}$ serves, not only as the objective function of the optimization problem, but as measure of the quality of the coupling as well. The numerical experiments presented in Figure \ref{fig:var_vs_q} is in perfect agreement with the hierarchy of inequalities (\ref{inequality_1qN}).

\subsection{Numerical results}\label{section:numerics_macro}

In this section we will discuss numerical results that compare the various couplings presented in the previous sections. As a reference model we use a combination of adsorption/desorption (\ref{ising_new_configuration}),(\ref{ising_rates}) and diffusion (\ref{diffusion_new_configuration}),(\ref{diffusion_rates}) mechanisms, see  (\ref{ising_diffusion_generator}). As observable we used the coverage defined by (\ref{coverage}). The parameters used in this model are $\beta=0.1, J=1, h=0, c_a=c_d=c_{diff}=1,N=100$ and the final time $T=10$. The perturbation is done in the $\beta$ parameter with $\epsilon=10^{-3}$. The estimator (\ref{estimator_mean}) is computed using $2000$ samples. 

As suggested by inequality (\ref{inequality_1qN}) in Section \ref{sec:macroscopic_coupling} the variance reduction using the macroscopic generator (\ref{macro_coupled}), corresponding to $q=N$, is $20$ times more than that obtained using the microscopic generator (\ref{generator_cpl}), corresponding to $q=1$, as seen in Figure \ref{fig:var_vs_q_vs_t}. On the other hand, simulating the coupled algorithm for $q=1$ and $q=N$ is $3.8$ and $2.4$ times slower, respectively, than simulating two uncoupled processes or the CRN method, see Table \ref{table:execution_time}.  Moreover, this simulation shows that the variance reduction is still obtained for large times even though the design of the algorithm is based on the assumption that $t\ll 1$, see discussion in Section \ref{sec:microscopic_coupling}.

In Figure \ref{fig:var_vs_q} a numerical verification of inequality (\ref{inequality_1qN}) is presented showing that the maximum variance reduction is obtained for $q=N$, at least for the coverage used as the observable in the calculation in Appendix \ref{appendix:inequality_proof}. The first point in the graph, $q=0$, corresponds to the variance ratio of the uncoupled process (\ref{simplest_coupling}) and is shown here only for comparison purposes. The overall variance reduction for the macroscopic coupling, $q=N$,  is $200$ times more than using two uncoupled processes or the CRN method. Thus, in order to get the same variance with the macroscopic coupling, the uncoupled algorithm needs about $80$ times more computational time than the macroscopically coupled process.
This result combined with the execution time comparisons of Table \ref{table:execution_time} shows that the optimal choice for $q$ is when $q=N$.
The results presented in this section can be reproduced by Matlab scripts found in\cite{matlab_code}, see also Appendix \ref{app:implementation}.

\begin{table}[htpb]
\begin{center}
\begin{tabular}{|c||c|c|c|c|c|c|c|c|c|}
\hline
$q$  &      1  &   2 &    4  &   5  &  10  &  20  &  25  &  50  & 100   \\ \hline 
$\frac{\textrm{coupled execution time}}{\textrm{uncoupled execution time}}$ &   3.8 &    3.1 &   2.8  &  2.7  &  2.5  &  2.4 &   2.4  &  2.4 &   2.4 \\ \hline
\end{tabular}
\end{center}
\caption{Ratio of coupled to uncoupled averaged execution time for the model presented in Section \ref{section:numerics_macro}.}
\label{table:execution_time}
\end{table}


\section{Observable-based coupling methods and gradient-free  information theoretic tools.}
The coupling sensitivity analysis methods  presented  in earlier sections focused on the
quantification of parameter sensitivities for specific observables of the KMC process.
However, in many 
complex spatial  reaction-diffusion  KMC
algorithms  there is a  combinatorial explosion in
the number of parameters \cite{ref3}. The high-dimensional parameter space creates  a seemingly intractable challenge for any finite-difference, {\em gradient-based}, sensitivity analysis method,
including the low-variance coupling methods we develop here, due to the fact that all partial derivatives (\ref{partial_derivative}) need to be computed.

On the other hand in  \cite{PK_RER} the authors developed a computationally tractable, {\em gradient-free} sensitive analysis method
suitable for complex stochastic dynamics, including spatial KMC algorithms, as well as systems with a very high-dimensional parameter space such as biochemical reaction networks with over $200$ parameters \cite{PKV}. The method is 
based on developing {\em computable}  information theory-based sensitivity metrics such as relative entropy, relative entropy rate  and Fisher Information Matrix  
at the path-space level, i.e., studying the sensitivity of  the entire stationary time-series.
Overall, in physicochemical and biological models, we focus  on specific 
observables such as coverage, populations, spatial  correlations, population variance, while important observables such as  autocorrelations and exit  times depend on 
 the entire time series.
%
%
%
%
%
%
%
%
%
%
%
%
%
%
Therefore, it is
plausible to attempt to connect the parameter sensitivities of  observables to the gradient-free, e.g. information-theoretic
methods proposed  \cite{PK_RER}. Indeed, relative entropy can provide an upper bound for a large family
of observable functions through the  Csiszar-Kullback-Pinsker inequality \cite{Cover:91}. 
More precisely, for any  bounded observable function $f$, the Pinsker inequality states that
\begin{equation}
|\mathbb E_{Q^{\theta}}[f] - \mathbb E_{Q^{\theta+\epsilon}}[f]| \leq \max_{\sigma}|f(\sigma)| \sqrt{2 \RELENT{Q^{\theta}}{Q^{\theta+\epsilon}}}\, ,
\label{Pinsker:ineq}
\end{equation}
and $\RELENT{Q^{\theta}}{Q^{\theta+\epsilon}}$ is the relative entropy between the {\em path space distributions}
$Q^{\theta}$ and $Q^{\theta+\epsilon}$ of the time series with parameters $\theta$ and $\theta+\epsilon$ respectively. Practically, the relative entropy $\RELENT{Q^{\theta}}{Q^{\theta+\epsilon}}$
quantifies the loss of information in the time-series distribution due to  a perturbation from $\theta$ to $\theta +\epsilon$. In \cite{PK_RER} it is shown that it is a computable observable and that it can be approximated by a corresponding Fisher Information Matrix (FIM)
defined on path space:
\begin{equation}
\RELENT{{Q^{\theta}}}{{Q^{\theta+\epsilon}}} = \frac{1}{2} \epsilon^T \FISHERR(Q^{\theta}) \epsilon + O(|\epsilon|^3)
\label{GFIM}
\end{equation}
where $\FISHERR(Q^{\theta})$ is a  $k\times k$ matrix--$k$ is the dimension of the parameter vector $\theta$--and  can be considered as a path-wise analogue for the classical Fisher Information
Matrix (FIM). 
Moreover, we have  an explicit formula for the path-wise FIM involving only the transition rates hence  it
is numerically computable as an observable of the process \cite{PK_RER}.
As in the classical FIM for parametrized distributions \cite{Cover:91}, \VIZ{GFIM} is the Hessian of the RER
which geometrically corresponds to the curvature around the minimum value of the
relative entropy. Therefore, the spectral analysis of the FIM yields a derivative-free sensitivity analysis method,
characterized by the the eigenvalues and eigenvectors of the FIM: the eigenvectors   correspond to the (hidden) sensitivity directions of the system, 
and  they are ordered  in terms of the eigenvalues from most  (higher eigenvalues) to least (lower eigenvalues) sensitive.

An outcome of \VIZ{Pinsker:ineq}  is that if the (pseudo-)distance between two 
distributions defined  by $\RELENT{Q^{\theta}}{Q^{\theta+\epsilon}}$ is controlled,
then the error between the two distributions is also controlled for any bounded
observable.
In the context of sensitivity analysis, inequality \VIZ{Pinsker:ineq} implies  that if the
relative entropy is small, i.e., insensitive in a particular parameter  direction, then,
any  bounded observable $f$ is also expected to be  insensitive towards the same
direction. 
%
%
Thus,   \VIZ{Pinsker:ineq}, combined with our current work on coupling  suggest  the following strategy  for sensitivity analysis:
(a) first,   the upper bound in \VIZ{Pinsker:ineq} constitutes a theoretical indicator that relative entropy is a reliable  tool for sensitivity analysis; (b)  
from a practical perspective, 
\VIZ{Pinsker:ineq} can   rule out insensitive directions in parameter space, given by eigenvectors of the FIM in \VIZ{GFIM};  in turn, this fact provides a  significant
advantage  in  the study of  models with a very large number of parameters, where the calculation of all gradients, even using coupling would be impractical; (c) finally, the remaining sensitive directions can be 
 explored using our proposed coupling methods.

%

\section{Conclusions}

In this paper we proposed a coupling method for variance reduction  of finite difference-based sensitivity analysis for  lattice kinetic Monte Carlo algorithms. Variance reduction in sensitivity analysis is a challenge of particular importance in high-dimensional stochastic dynamics, such as the ones addressed here, since the high computational cost of individual realizations of the stochastic process renders prohibitive the generation of a large number of samples for reliable  ensemble averaging. Our proposed method relies  on defining a new coupled continuous time Markov Chain through a suitable generator that acts on observables of the involved high-dimensional stochastic processes associated with parameters $\theta$ and $\theta+\epsilon$. 
The rates of the coupled generator--and the corresponding coupled stochastic processes--are obtained by solving an optimization problem associated with  minimizing  the  variance between the coupled stochastic dynamics. Moreover,  
the optimization problem depends on the particular  observable quantity of interest. Thus, the form of the coupled rates depends directly on the choice of the observable function, hence  it is a ``goal-oriented'' method.  
The implementation of the coupled method is a BKL--type algorithm  in the sense that events are categorized into pre-defined sets. In the classical BKL algorithm \cite{BKL} the events are divided into classes  of equal rates, i.e. level sets of the rates. However here the events are divided into classes depending on the targeted observable's  level sets, since we are interested in tight couplings of time series of specific observables. 
The division of the state space into classes is being done automatically by the provided code, overcoming the problem of predefining the classes.
Numerical examples of spatial KMC, e.g. adsorption/desorption/diffusion processes, presented throughout  the paper show that the variance reduction can be improved by  two orders of magnitude compared to typical  methods used up to now, such as the Common Random Number method. At the same time the computational overhead of constructing the coupled rates is two times slower than that of the CRN method leading to an overall speed up factor of two orders of magnitude. Furthermore, numerical experiments also demonstrate that the optimization functional (\ref{optimization_functional_2}) indeed constitutes a diagnostic tool for the design and evaluation of different couplings, e.g., see (\ref{inequality_1qN}) and Figure \ref{fig:var_vs_q}. 
Finally, we noted that even  the proposed coupling methods  for sensitivity analysis are not an  efficient approach for systems with a very high dimensional parameter space, due to the fact that coupling  is a gradient method requiring the calculation of all (discrete) partial derivatives. Moreover, the results presented in this paper are all reproducible by Matlab scripts\cite{matlab_code}. We expect that a  combination of the proposed low variance coupling methods with gradient-free methods, such as the Relative Entropy Rate method \cite{PK_RER}, see \VIZ{Pinsker:ineq},  can also provide  a realistic and systematic approach  towards the  sensitivity analysis of  systems with a high-dimensional parameter space.

\section*{Acknowledgements}
The work of GA was supported by the Office of
Advanced Scientific Computing Research, U.S. Department
of Energy under Contract No. DE-SC0002339 and by   the European Union
(European Social Fund) and Greece (National Strategic Reference Framework),
under the THALES Program, grant AMOSICSS.
The work of MAK    was supported in part by the Office of
Advanced Scientific Computing Research, U.S. Department
of Energy under Contracts No. DE-SC0002339 and No. DE-SC0010723.



\appendix

\section{Examples of lattice models}\label{app:examples}

In order to lighten the notation in the following examples, let us give a helpful definition.
\begin{definition}
Let $A$ be a set. Then,
\begin{equation}\label{set_notation}
A[\textbf{condition}] = 
 \begin{cases}
A, & \textrm{if}  \textbf{ condition} \textrm{ is satisfied}, \\
\emptyset, & \textrm{else}
  \end{cases} \PERIOD
\end{equation}
\end{definition}

\noindent
\textbf{\uline{Example 1}} The ZGB model \cite{ZGB} is a simplified atomistic lattice-gas model of CO oxidation and consists of three
reactions,
\begin{equation*}
 \begin{split}
  \textrm{CO} &\textrm{ adsorption on an empty site} \\
  \mathrm{O}_2 &\textrm{ adsorption on two neigbor empty sites} \\
  \mathrm{CO} + \mathrm{O} &\longrightarrow \mathrm{ CO}_2  \textrm{ and instantaneous desorption leaving two empty neigbor sites}\PERIOD
 \end{split} 
\end{equation*}
In this case $\Sigma=\{-1,0,1\}$, where $-1,0$ and $1$ correspond to CO, vacant site and O respectively. 
The neighborhood $\Omega_x$ is described in Figure \ref{fig:neighborhood}. Let us  also  define $\sigma_i:=\sigma(x_i)$. Then the possible transitions are given through the rate function,
\begin{equation}\label{ZGB_rates}
 c(x;\sigma,\omega) = 
 \begin{cases}
  c_a, & \textrm{if } \omega\in\Omega_1(x,\sigma), \\
 1- c_a, & \textrm{if } \omega\in\Omega_2(x,\sigma), \\
  c_r, & \textrm{if } \omega\in\Omega_3(x,\sigma), \\
    0,   & \textrm{else}\COMMA
 \end{cases} 
\end{equation}
where using notation (\ref{set_notation}) we write,
\begin{equation*}
\begin{split}
 \Omega_1(x,\sigma) &= \Bigl\{(-1,\sigma_2,\ldots,\sigma_5)\Bigr\} \Bigl[\sigma_1=0 \Bigr] \\
 \Omega_2(x,\sigma) &= \bigcup_{i=1}^4  \Bigl\{ (1,\sigma_2,\ldots,\omega_i,\ldots,\sigma_5) \textrm{ with } \omega_i=1 \Bigr\} \Bigl[\sigma_1=0 \textrm{ and } \sigma_i=0 \Bigr] \\ 
 \Omega_3(x,\sigma) &=\bigcup_{i=1}^4  \Bigl\{ (0,\sigma_2,\ldots,\omega_i,\ldots,\sigma_5) \textrm{ with } \omega_i=0 \Bigr\}     \Bigl[   \sigma_1=\pm 1 \textrm{ and } \sigma_i=\mp 1  \Bigr]  \PERIOD
\end{split}
\end{equation*}

\noindent
\textbf{\uline{Example 2}} A more refined model of CO oxidation is presented in \cite{Evans}. This model has the same reactions with the
ZGB model with the addition of CO diffusion. In order to take into  account the strong O-O repulsions, the $\textrm{O}_2$ desorption is done
at diagonal nearest-neighbor sites with the additional constraint that the six sites adjacent to these are empty. The neighborhood of
the model is presented in Figure \ref{fig:neighborhood}.

In order to describe the sets of possible local configurations let us first define the set 
\begin{equation*}
 D_i=\{ \textrm{the indices of the adjacent neighbors at site } x_i\in\Omega_x \}, \; i=3,5,7,9 \COMMA
\end{equation*}
e.g.\ $D_3 =\{ 2,4,11,12\}$, see Figure \ref{fig:neighborhood}. Then the possible transitions are given through the rate function,
\begin{equation}\label{Evans_rates}
 c(x;\sigma,\omega) = 
 \begin{cases}
  c_a, & \textrm{if } \omega\in\Omega_1(x,\sigma), \\
 1- c_a, & \textrm{if } \omega\in\Omega_2(x,\sigma), \\
  c_{diff}, & \textrm{if } \omega\in\Omega_3(x,\sigma), \\
  c_d, & \textrm{if } \omega\in\Omega_4(x,\sigma), \\
  c_r, & \textrm{if } \omega\in\Omega_5(x,\sigma), \\
    0,   & \textrm{else}\COMMA
 \end{cases} 
\end{equation}
where using notation (\ref{set_notation}) we write,
\begin{equation*}
\begin{split}
\Omega_1(x,\sigma) &=  \Bigl\{(-1,\sigma_2,\ldots,\sigma_{17})\Bigr\}   \Bigl[\sigma_1=0 \Bigr] \\
\Omega_2(x,\sigma) &=\bigcup_{i=3,5,7,9} \Bigl\{ (1,\sigma_2,\ldots,\omega_i,\ldots,\sigma_{17}) \textrm{ with } \omega_i=1 \Bigr\}\cdots \\
 &\qquad\qquad\qquad\qquad\qquad \cdots\Bigl[ \sigma_1 =0 \textrm{ and } \sigma_i =0  \textrm{ and } \sigma_j=0,\; \forall j \in D_i  \Bigr ]  \\
\Omega_3(x,\sigma)  &= \bigcup_{i=2,4,6,8}  \Bigl\{ (0,\sigma_2,\ldots,\omega_i,\ldots,\sigma_{17}) \textrm{ with }\omega_i=-1 \Bigr\}
\Bigl[  \sigma_1=-1 \textrm{ and }\sigma_i=0  \Bigr]  \\
\Omega_4(x,\sigma)  &= \Bigl\{ (0,\sigma_2,\ldots,\sigma_{17}) \Bigr\}  \Bigl[\sigma_1=-1 \Bigr]  \\
\Omega_5(x,\sigma) &= \bigcup_{i=2,4,6,8}  \Bigl\{ (0,\sigma_2,\ldots,\omega_i,\ldots,\sigma_{17}) \textrm{ with } \omega_i=0 \Bigr\}  
\Bigl[  \sigma_1=\pm 1 \textrm{ and } \sigma_i=\mp 1  \Bigr]   \PERIOD
\end{split}
\end{equation*}
%

%

\section{Proof of inequality (\ref{inequality_1qN})}\label{appendix:inequality_proof}

In order to prove inequality  (\ref{inequality_1N}) we use the coverage (\ref{coverage}) as an observable and the sets (\ref{coverage_sets}) as the partition sets.
Thus the values of the $f(\sigma^{x,\omega}) - f(\sigma)$ in (\ref{optimization_functional_2}) are $-1,0,1$ for $k=1,2,3$ respectively. 
Using the rates (\ref{coarse_coupled_rates}) into the objective function (\ref{optimization_functional_2}) of the optimization problem (\ref{max_problem_general}), we have:
\begin{equation}\label{app_proof_1}
\begin{split}
\mathcal{F}[c_q;f] &=  
\sum_{j=1}^M\sum_{j'=1}^M \sum_{k=1,3} \sum_{k'=1,3} \sum_{(x,\omega)\in S_{k,j}(\sigma)} \sum_{(y,\omega')\in S_{k',j'}(\eta)} \ldots  \\
&\qquad \ldots \min \{ \lambda^A_{k,j}(\sigma) , \lambda^B_{k',j'}(\eta) \} \frac{c_A(x,\omega;\sigma)}{\lambda_{k,j}^A(\sigma)}  \frac{c_B(y,\omega';\eta)}{\lambda_{k',j'}^B(\eta)} \delta(k-k') \delta(j-j') \\
&=  \sum_{k=1,3}  \sum_{j=1}^M  \min \{ \lambda^A_{k,j}(\sigma) , \lambda^B_{k,j}(\eta) \}  \Bigl [ \sum_{(x,\omega)\in S_{k,j}} \sum_{(y,\omega')\in S_{k,j}}   \frac{c_A(x,\omega;\sigma)}{\lambda_{k,i}^A(\sigma)}  \frac{c_B(y,\omega';\eta)}{\lambda_{k,i}^B(\eta)} \Bigr ] \\
& = \sum_{k=1,3}  \sum_{j=1}^M  \min \{ \lambda^A_{k,j}(\sigma) , \lambda^B_{k,j}(\eta) \}   \COMMA
\end{split}
\end{equation}
where $\delta(x)=1$ if $x=0$ and zero otherwise.
The absence of $k=2$ is justified by the fact that $[f(\sigma^{x,\omega})-f(\sigma)][f(\sigma^{y, \omega'})-f(\sigma)]$  is zero for $(x,\omega)\in S_{k,j}(\sigma)$ and $(y,\omega')\in S_{k,j}(\eta)$  while when $k=1,3$ the same product equals to $1$.
By fixing $k$ in the summation of the last part of (\ref{app_proof_1}) we have: 
\begin{equation}\label{app_proof_2}
\begin{split}
\sum_{j=1}^M  \min \{ \lambda^A_{k,j}(\sigma)  &, \lambda^B_{k,j}(\eta) \} = \\
&= \sum_{j=1}^M \min \{ \sum_{(x,\omega)\in S_{k,j}(\sigma)} c_A(x,\omega;\sigma) , \sum_{(y,\omega')\in S_{k,j}(\eta)}  c_B(y,\omega';\eta) \} \\
&\leq  \min \{ \sum_{j=1}^M \sum_{(x,\omega)\in S_{k,j}(\sigma)} c_A(x,\omega;\sigma) ,  \sum_{j=1}^M \sum_{(y,\omega')\in S_{k,j}(\eta)}  c_B(y,\omega';\eta) \} \\
& = \min \{ \sum_{(x,\omega)\in S_k(\sigma)} c_A(x,\omega;\sigma) ,   \sum_{(y,\omega')\in S_k(\eta)}  c_B(y,\omega';\eta) \} \\
& =  \min \{ \lambda^A_{k}(\sigma) , \lambda^B_{k}(\eta) \} \PERIOD
\end{split}
\end{equation}
Finally, by summation over $k$ in both sides of (\ref{app_proof_2}) and that fact that 
\begin{equation}
\mathcal{F}[c_N;f]  = \sum_{k=1,3}  \min \{ \lambda^A_{k}(\sigma) , \lambda^B_{k}(\eta) \} \COMMA
\end{equation}
we get that
\begin{equation}
\mathcal{F}[c_q;f] \leq \mathcal{F}[c_N;f], \quad 1\leq q <N.
\end{equation}


\section{Implementation}\label{app:implementation}

The goal-oriented sensitivity analysis method proposed in this work was implemented in Matlab and the source code can be 
downloaded from here \cite{matlab_code}.
The code was written in a way such that  the user can easily specify and run different models: the neighborhood (see definition (\ref{def:neigh})), the set of all accessible local configurations (see definition (\ref{def:local_conf_omega})) and the corresponding rates (see also Appendix \ref{app:examples}) are described in two files with a specific pre-defined format. In a different file the preferred observable function can be defined by the user.  The microscopic (\ref{micro_rates_1}),(\ref{micro_rate_2}), macroscopic (\ref{macro_coupled_rates}) and all intermediate couplings (\ref{coarse_coupled_rates}) can be executed by changing the value of a variable.
Moreover, the code can easily handle 1D as well as 2D models.

The code is split into multiple functions so that  new methods can be easily integrated. For example, the current version of the code includes only the the Stochastic Simulation Algorithm \cite{Gillespie76} but a different method for the time advancement of the stochastic system can be included.
Finally, all examples and figures presented in Sections \ref{section:numerics_micro} and \ref{section:numerics_macro} can be reproduced by the included Matlab scripts \cite{matlab_code}.


\end{document}